\documentclass[ijoc,nonblindrev]{oo}

\OneAndAHalfSpacedXII % Current default line spacing

\usepackage{natbib}
\bibpunct[, ]{(}{)}{,}{a}{}{,}%
\def\bibfont{\small}%

\TheoremsNumberedThrough     % Preferred (Theorem 1, Lemma 1, Theorem 2)

\ECRepeatTheorems

\EquationsNumberedThrough    % Default: (1), (2), ...

\RUNAUTHOR{Daryalal and Pouya}

\RUNTITLE{Network Migration Problem: An LBBD Approach Driven by CG and CP}

\TITLE{Network Migration Problem: A Logic-based Benders Decomposition Approach Driven by Column Generation and Constraint Programming}

\ARTICLEAUTHORS{%
\AUTHOR{Maryam Daryalal$^\star$, Hamed Pouya$^\dagger$}
 \AFF{Department of Mechanical and Industrial Engineering, University of Toronto, Toronto, Ontario M5S 3G8, Canada 
 }
 \AFF{$^\star$\EMAIL{daryalal@mie.utoronto.ca}, $^\dagger$\EMAIL{h.pouya@utoronto.ca}
 }
}

\ABSTRACT{%
Telecommunication networks frequently face technological advancements and need to upgrade their infrastructure. Adapting legacy networks to the latest technology requires synchronized technicians responsible for migrating the equipment. The goal of the network migration problem is to find an optimal plan for this process. This is a defining step in the customer acquisition of telecommunications service suppliers, and its outcome directly impacts the network owners’ purchasing behaviour. We propose the first exact method for the network migration problem, a logic-based Benders decomposition approach that benefits from a hybrid constraint programming-based column generation in its master problem and a constraint programming model in its subproblem. This integrated solution technique is applicable to any integer programming problem with similar structure, most notably the vehicle routing problem with node synchronization constraints. Comprehensive evaluation of our method over instances based on six real networks demonstrates the computational eﬀiciency of the algorithm in obtaining quality solutions. We also show the merit of each incorporated optimization paradigm in achieving this performance.}


\usepackage{amsmath,amssymb}
\usepackage{color}
\usepackage{graphicx}
\usepackage{subfig}
\usepackage{url}
\usepackage{multirow}
\usepackage{array}
\usepackage{amsfonts}
\usepackage{mathtools}

\usepackage{enumitem}
\usepackage{lineno}
\usepackage{tikz}
\usepackage{algorithm}
\usepackage[noend]{algpseudocode}
\usepackage{booktabs}
\usepackage{longtable}
\usepackage{geometry}
\usepackage{xr-hyper}
\usepackage[colorlinks=true, linkcolor=blue, citecolor=blue]{hyperref}

\def \checkmark{\tikz\fill[scale=0.4](0,.35) -- (.25,0) -- (1,.7) -- (.25,.15) -- cycle;} 
\modulolinenumbers[5]

\def \durshift {\Delta_{\textsc{shift}}}
\def \g {\gamma}
\def \G {\Gamma}
\def \Gp {\Gamma'}
\def \duration {\Delta}
\def \durmigr {\theta}
\def \mwDur {\underline{\duration}}
\def \travel {T}
\def \dual {\pi}
\def \aEng {\alpha^{\textsc{eng}}}
\def \cirSS {\phi_{ss'}}
\def \maxTechR {\eta_r^{\textsc{tech}}}
\def \maxCir {\eta^{\textsc{cir}}}
\def \maxEng {\eta^{\textsc{eng}}}
\def \siteSet {\mathcal{S}}
\def \regSet {\mathcal{R}}
\def \cirSet {\mathcal{C}}

\def \costT {\textsc{cost}^{\textsc{tech}}}
\def \costE {\textsc{cost}^{\textsc{eng}}}
\def \cost {\textsc{cost}}
\def \sitePair {\mathcal{S}_p}
\def \src {\textsc{src}}
\def \dst {\textsc{dst}}
\def \SRC {\textsc{src}}
\def \DST {\textsc{dst}}
\def \wSet {\mathcal{W}}
\def \traTime {T}
\def \T {\mathcal{T}}

\def \BDMP {\text{MP}^{\textsc{BD}}}
\def \BDSP {\text{SP}^{\text{BD}}_w(\overline{\boldsymbol{m}}_w)}
\def \BDSPSol {\overline{\text{SP}}^{\text{BD}}_w(\overline{\boldsymbol{m}}_w)}
\def \LBBDSP {\text{SP}^{\text{LBBD}}_w(\overline{\boldsymbol{m}}_w)}
\def \CGMP {\text{RMP}^{\text{CG}}_w}
\def \CGSP {\text{SP}^{\text{CG}}_{rw}}
\def \CGSPSol {\overline{\text{SP}}^{\text{CG}}_{rw}}
\def \CP {\text{CP}_w(\overline{\boldsymbol{m}}_w)}
\def \CPSol {\overline{\text{CP}}_w(\overline{\boldsymbol{m}}_w)}
\def \LBBDMP {\text{MP}^{\textsc{LBBD}}}
\def \iter {\tau}
\def \greater {\kappa}
\def \optcuts {\Lambda^{\textsc{opt}}_w}
\def \feascuts {\Lambda^{\textsc{feas}}_w}
\def \artvar {\rho}

\def \techsiteTS {x_{rts}}
\def \seqT {\textsc{seq}_{rt}}
\def \numTS {n^{\textsc{s}}_{rts}}
\def \numTSSp {n^{\textsc{sp}}_{rtss'}}
\def \numTSpS {n^{\textsc{sp}}_{rts's}}
\def \worktimeT {\textsc{wtime}_{rt}}
\def \techShiftDurT {\Delta_{\textsc{shift}}^{rt}}

\def \IntervalVar {\text{IntervalVar}}
\def \SequenceVar {\text{SequenceVar}}
\def \Span {\text{Span}}
\def \NoOverlap {\text{NoOverlap}}
\def \IntegerVar {\text{IntegerVar}}
\def \IfThen {\text{IfThen}}
\def \PresenceOf {\text{PresenceOf}}
\def \EndOf {\text{EndOf}}
\def \LengthOf {\text{LengthOf}}

\renewcommand{\proof}[1]{\textit{Proof}. #1\hfill\Halmos}
\newcommand{\proofNoH}[1]{\textit{Proof}. #1}

\newcolumntype{C}[1]{>{\centering\arraybackslash}p{#1}}

\def \EUnet {\textsc{EUNetworks}}
\def \Sago {\textsc{Sago}}
\def \Savvis {\textsc{Savvis}}
\def \Vision {\textsc{VisionNet}}
\def \VisionNet {\Vision}
\def \NextGen {\textsc{NextGen}}
\def \Pionier {\textsc{Pionier}}

\makeatletter
\def\LT@start{%
  \let\LT@start\endgraf
  \endgraf\penalty\z@\vskip\LTpre\endgraf
  \dimen@\pagetotal
  \advance\dimen@ \ht\ifvoid\LT@firsthead\LT@head\else\LT@firsthead\fi
  \advance\dimen@ \dp\ifvoid\LT@firsthead\LT@head\else\LT@firsthead\fi
  \advance\dimen@ \ht\LT@foot
  \dimen@ii\vfuzz
  \vfuzz\maxdimen
    \setbox\tw@\copy\z@
    \setbox\tw@\vsplit\tw@ to \ht\@arstrutbox
    \setbox\tw@\vbox{\unvbox\tw@}%
  \vfuzz\dimen@ii
  \advance\dimen@ \ht
        \ifdim\ht\@arstrutbox>\ht\tw@\@arstrutbox\else\tw@\fi
  \advance\dimen@\dp
        \ifdim\dp\@arstrutbox>\dp\tw@\@arstrutbox\else\tw@\fi
  \advance\dimen@ -\pagegoal
  \ifdim \dimen@>\z@\vfil\break\fi
      \global\@colroom\@colht
  \ifvoid\LT@foot\else
     \global\advance\vsize-\ht\LT@foot
    \global\advance\@colroom-\ht\LT@foot
    \dimen@\pagegoal\advance\dimen@-\ht\LT@foot\pagegoal\dimen@
    \maxdepth\z@
  \fi
  \ifvoid\LT@firsthead\copy\LT@head\else\box\LT@firsthead\fi\nobreak
  \output{\LT@output}}

\def\endlongtable{%
  \crcr
  \noalign{%
    \let\LT@entry\LT@entry@chop
    \xdef\LT@save@row{\LT@save@row}}%
  \LT@echunk
  \LT@start
  \unvbox\z@
  \LT@get@widths
  \if@filesw
    {\let\LT@entry\LT@entry@write\immediate\write\@auxout{%
      \gdef\expandafter\noexpand
        \csname LT@\romannumeral\c@LT@tables\endcsname
          {\LT@save@row}}}%
  \fi
  \ifx\LT@save@row\LT@@save@row
  \else
    \LT@warn{Column \@width s have changed\MessageBreak
             in table \thetable}%
    \LT@final@warn
  \fi
  \endgraf\penalty -\LT@end@pen
  \ifvoid\LT@foot\else
    \global\advance\vsize\ht\LT@foot
    \global\advance\@colroom\ht\LT@foot
    \dimen@\pagegoal\advance\dimen@\ht\LT@foot\pagegoal\dimen@
  \fi
  \endgroup
  \global\@mparbottom\z@
  \endgraf\penalty\z@\addvspace\LTpost
  \ifvoid\footins\else\insert\footins{}\fi}

\def\LT@output{%
  \ifnum\outputpenalty <-\@Mi
    \ifnum\outputpenalty > -\LT@end@pen
      \LT@err{floats and marginpars not allowed in a longtable}\@ehc
    \else
      \setbox\z@\vbox{\unvbox\@cclv}%
      \ifdim \ht\LT@lastfoot>\ht\LT@foot
        \dimen@\pagegoal
\advance\dimen@\ht\LT@foot
        \advance\dimen@-\ht\LT@lastfoot
        \ifdim\dimen@<\ht\z@
          \setbox\@cclv\vbox{\unvbox\z@\copy\LT@foot\vss}%
          \@makecol
          \@outputpage
            \global\vsize\@colroom
          \setbox\z@\vbox{\box\LT@head}%
        \fi
      \fi
      \unvbox\z@\ifvoid\LT@lastfoot\copy\LT@foot\else\box\LT@lastfoot\fi
    \fi
  \else
    \setbox\@cclv\vbox{\unvbox\@cclv\copy\LT@foot\vss}%
    \@makecol
    \@outputpage
      \global\vsize\@colroom
    \copy\LT@head\nobreak
  \fi}
\makeatother

\begin{document}

\maketitle

% ---------------------------------------------------------
% ---------------------------------------------------------
% ---------------------------------------------------------	

\section{Introduction}
\label{sec:introduction}
%HP Motivation
In telecommunication industries, network migration is the process of upgrading the existing infrastructure of a deployed network. A telecommunication network  is composed of a set of sites (demand points), and circuits that transmit the traffic between the sites.  Migration of such a network is performed by upgrading the circuits one by one. In order to upgrade every circuit, two synchronized technicians migrate its two endpoints within a time window.  The goal of the \emph{network migration problem} (NMP) is to find the upgrade order of these circuits such that the associated costs are minimized. 
Migration of a network is  a strategic decision that can lead to immense savings of  10 to 100 times in power and space as reported by \cite{ciena2013}. Nevertheless,  the process is quite costly and complex, with some circuits stretching over a continent. Furthermore, every circuit migration comes with a disruption as the endpoints disconnect from the equipment, thus affecting the efficiency of the network and the customers' satisfaction of the migration solution, i.e., the frequency of disruptions they experience during the upgrade. Consequently, a well-crafted plan is critical for the success of the migration. 
In this paper, we present the first exact solution method for the NMP, a \emph{logic-based benders decomposition} (LBBD) algorithm that integrates Benders Decomposition (BD), Column Generation (CG) and Constraint Programming (CP),  enabling us to decompose the decision space into smaller subsets, each amenable to one of these solution frameworks.

%Existing papers on Network Migration Problem
The NMP can also be stated in the context of the \emph{vehicle routing problem with synchronized constraints} (VRPS).
In the VRPS, at least one vertex or arc requires simultaneous visits of vehicles, or successive visits resulting from some precedence constraints \citep{eksioglu2009vehicle}. The synchronization constraints can be over the arcs (the synchronized arc routing problem, SARP), or the nodes (the synchronized node routing problem, SNRP).  
%NMP is SNRP
In the NMP, vehicles and nodes correspond to technicians and sites, respectively.
Since  the tasks (i.e., upgrading the circuit endpoints) are defined over the nodes, the NMP is a special case of the SNRP and an application of the VRPS in the telecommunications domain.
% SNRP
In the SNRP, typically some nodes need to be visited by more than one vehicle because the personnel do not have the same expertise. 
% with instances of up to 45 customers and 16 vehicles. 
Some examples  are: \cite{bredstrom2007abranch} for a vehicle routing and scheduling problem,  \cite{reinhardt2013synchronized} for the airport transportation,  \cite{labadie2014iterated} and \cite{hashemi2020vehicle} for home healthcare systems, \cite{hojabri2018large} for an SNRP with precedence constraints and synchronization of two types of vehicles, and \cite{li2020branch} for a variation of the SNRP where customers have multiple options for time windows.
There are also studies considering both the SARP and the SNRP at the same time.
\cite{salazar2013synchronized} studied the road marking operations as a synchronized arc and node routing problem such that several capacitated vehicles are used to paint the lines on the roads and a tank vehicle is used to replenish the painting vehicles. Interested reader may refer to \cite{drexl2012synchronization} for a review on VRPS problems and their classifications. 
In regards to the applicability of these studies to the NMP, an important issue is their \emph{local} synchronisation assumption, meaning that the arrival of vehicles is synchronized at the same node.
    However, the NMP involves several sets of technicians distributed over multiple regions
    which can be synchronized with several other technicians from the same or other regions based on the location of the circuit endpoints. Besides, \cite{li2020branch} is the only exact method that considers multiple depots and time windows (the same as the NMP), but the size of the instances it can solve is up to 40 customers with 3 time windows which is far from the need of telecommunication networks.
Therefore, the existing works in the literature of the VRPS are not suited for the NMP's level of complexity. 

As for the telecommunications literature, there are different cost models and strategies on the suitable time and technology to migrate a telecommunication network \citep{podhradsky2004, almughaless2010optimum, turk2012networkMigOpt,poularakis2019optimizing}. Moreover, study of the operational aspects of the NMP has recently gained momentum. \cite{ble13} studied the migration of a network as the problem of finding the order of link upgrades that minimizes the total disruption time. 
With no travels allowed between the sites, their problem involves no routing decisions.
\cite{brigitteHamed} proposed a CG-based heuristic for the NMP that decides on the technician-to-circuit assignments, as well as travel paths for the technicians. \cite{pouya2017efficient}  observed that the prior mathematical formulations suffer from the highly symmetric nature of the circuits and technicians. A new symmetry-breaking model  was developed, leading to significant gains in computational effort. Yet again there is no guarantee that the obtained solutions are indeed optimal/feasible for the individual circuit assignments.  Subsequently,  \cite{jaumard2018migration} designed a two-phase CG-based 
heuristic for the planning of real-size networks. A greedy algorithm was proposed by  \cite{dawadi2021migration} for the cost minimization of a multi-technology migration problem. With the objective of minimizing the number of out-of-service sites and travels, \cite{javad2021efficient} modeled the NMP as a binary quadratic program. They derived heuristic solutions using Digital Annealer and solved instances with up to 64 circuits.  

% Note that, a ``plan'' does not have the exact schedule of the technicians, rather it makes sure that both endpoints of a circuit are migrated within a given time window. 
% Planning helps the network operators have an insight about the minimum migration costs. 

In this work, we develop an exact decomposition framework for the NMP. To the best of our knowledge, both in the context of the telecommunications problems and as a VRPS,  this is the first method that exactly solves the NMP with a certificate of optimality/infeasibility. We decompose the NMP into three problems and link them all by designing an LBBD algorithm.  By doing so, we are able to leverage the power of different solution techniques for linear programming and combinatorial optimization, and delegate the task of solving each problem to the most suitable optimization paradigm.  Numerical experiments defined over various real networks demonstrate the effectiveness of our algorithm in obtaining quality solutions with reasonable computational effort. Given that the NMP can be viewed as an SNRP, our method can also be adapted to a wide class of integer programming problems.

\noindent\textbf{Contributions}. The contributions of our work are summarized as follows.
\begin{itemize}
    \item We develop the first exact solution method for the network migration problem in order to find the optimal planning solutions, i.e.,  the order of circuit upgrades, along with the technician assignment and routing decisions. Although our solution framework is developed for the NMP, it is also applicable for the VRP with node synchronization.
    \item In order to reduce the computational effort of the LBBD master problem, we further decompose it via the Benders decomposition, resulting in a mixed-integer program (MIP) as the master problem, and a CP-based CG formulation as its subproblem. Additionally, we augment this CG formulation with an auxiliary MIP subproblem, leading to a hybrid CP/MIP-based CG model that significantly improves its performance.
    \item Considering the planning nature of our LBBD subproblem, we propose a CP model that, given the number of migrated circuits in a maintenance window, decides on the optimal technician assignments, order of circuit upgrades, and travel routes, if any.
    \item For our LBBD-based decomposition framework, we design valid feasibility and optimality cuts that guarantee the convergence of the algorithm. We also characterize a set of solutions other than the candidate for which our optimality cuts are tight.
    \item We evaluate the proposed LBBD algorithm on instances defined over six real backbone and regional networks and provide detailed algorithmic analysis and discussions on the implementation choices, along with managerial insights on the trade-offs among the migration cost, resource usage, and the duration of the migration. 
\end{itemize}
The remainder of the paper is organized as follows.
In Section \ref{sec:prob} we present the problem statement and its mathematical formulation. In Section \ref{sec:solution-method} we develop an LBBD solution framework for solving the NMP. In
Section \ref{sec:numer_results} we evaluate the performance of our algorithm on benchmark networks and provide  managerial insights. Section \ref{sec:conclusion} concludes the paper.

% ---------------------------------------------------------
% ---------------------------------------------------------
% ---------------------------------------------------------

\section{Problem Description}\label{sec:prob}
In this section, we formally describe the problem and introduce the sets and parameters. Next, we formulate the NMP as a CG-based integer linear program (ILP).
\subsection{Problem Statement}
\begin{figure}[t]
    	\centering
    	\includegraphics[scale = 0.11]{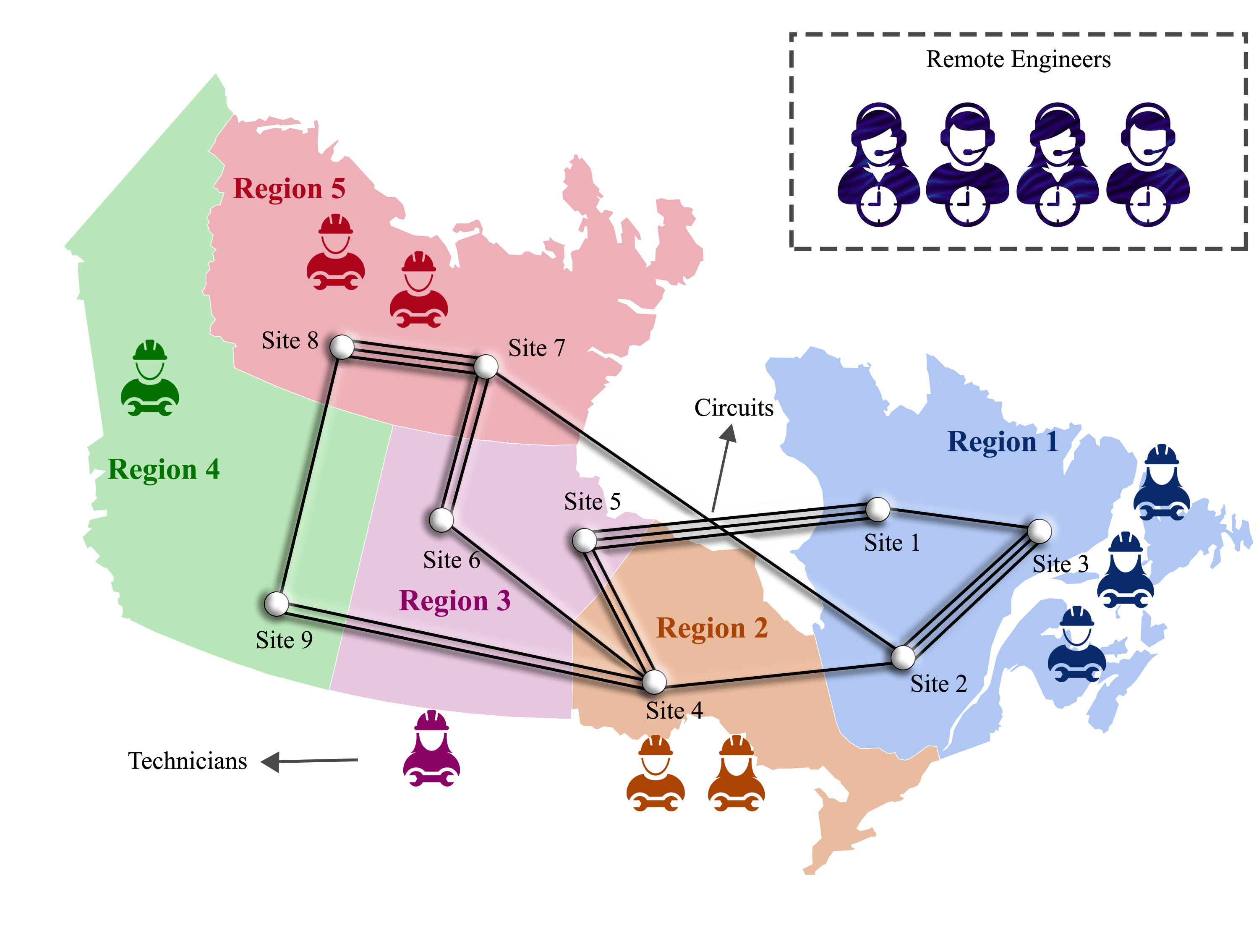}
    	\caption{An example of a telecommunication network for the network migration problem}
    	\label{fig:statement}
    \end{figure}

    The network migration problem is defined on a telecommunication network represented by a set of \emph{sites} $\siteSet$ and a set of \emph{circuits} $\cirSet$ between the site pairs $ \{s,s'\}\in\sitePair $.  Every site $ s $ is located in a geographical \emph{region} $r\in \regSet$ (e.g., a city). To each region $r$, $\maxTechR$ number of \emph{technicians} are assigned. Technicians assigned to a given region $r$ can only work in that region. Considering that  migrations often occur during low traffic periods (mostly nights), and the distance between the regions may require long-haul flights, employing local technicians is the safest option to minimize unforeseen impacts of the travels on migration planning. A circuit is migrated by disconnecting its endpoints from the old equipment and connecting them to the new one. These operations are performed by two technicians within the same maintenance window, each working at one endpoint. Additionally, every circuit migration requires an \emph{engineer} that coordinates the technicians remotely and does not need to be present in the working site. There are at most $\maxEng$ engineers available, and every engineer can coordinate up to $\aEng$ technicians.  Figure \ref{fig:statement} demonstrates an example of such a network with 5 regions, 9 sites, and 20 circuits.
    Migration of the network is performed during a \emph{maintenance window}, which is a period of time usually at night or a specially low traffic time on the network and also the roads (in case of travel between the sites). Every maintenance window $w\in \wSet$ has a predefined duration e.g., 8 hours, and all operations have to be completed within this duration. Time required to migrate a circuit $c\in \cirSet$ is $\durmigr$. Since migrating every circuit results in a short disruption in the network and the number of disruptions cannot violate clients' Service Level Agreement (SLA) \citep{fawaz2004service}, there is a limit $\maxCir$ on the number of migrated circuits per maintenance window.

\begin{figure}[t]
    	\centering
    	\includegraphics[scale = 0.11]{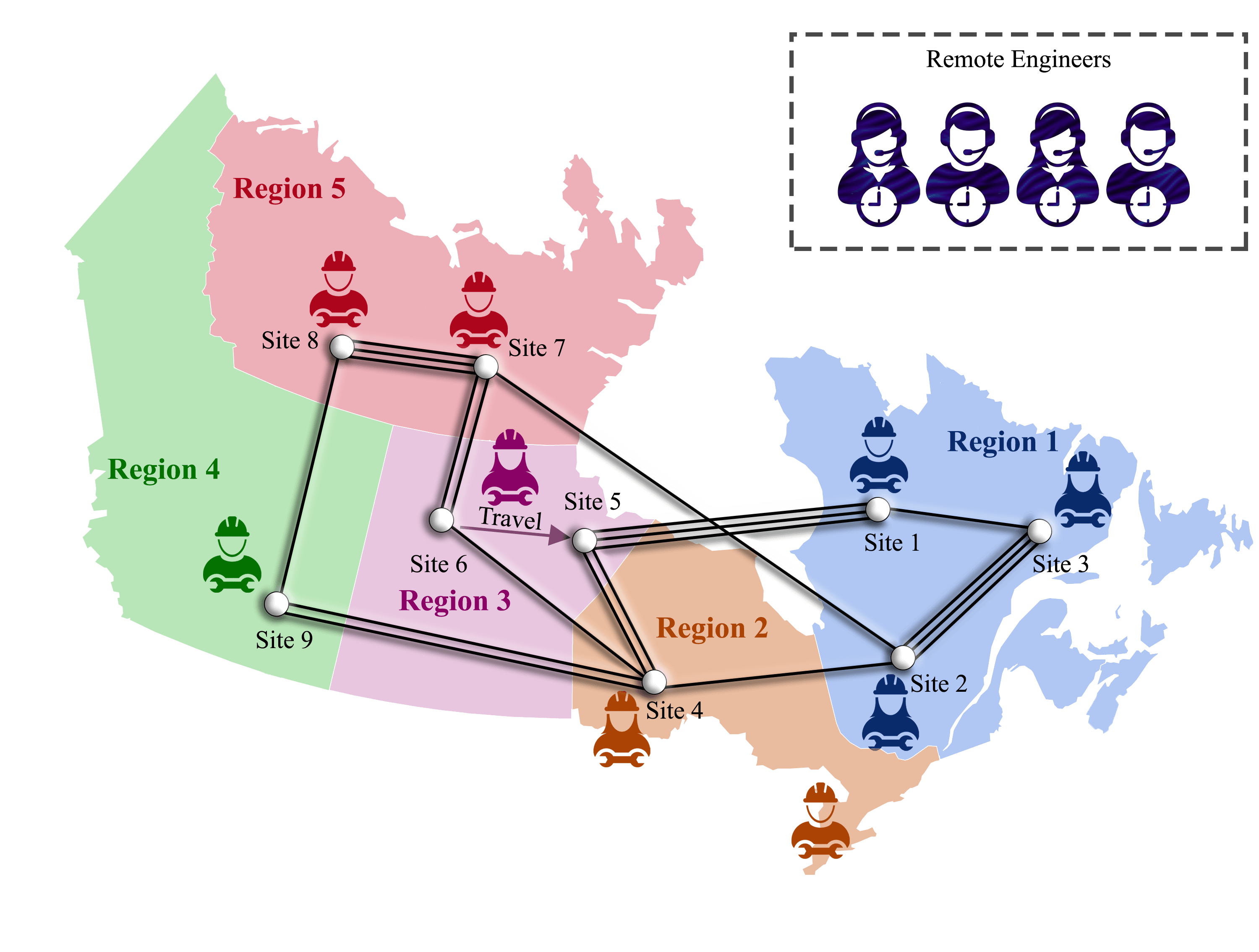}
    	\caption{A possible solution as a subset of shifts}
    	\label{fig:configuration}
    \end{figure}

    A technician working in a given region $r$ during maintenance window $w$ is responsible for a \emph{shift}. A shift is defined as a set of circuit endpoints migrated by a single technician during a maintenance window, together with any travels between the sites. Figure \ref{fig:configuration} represents a subset of possible shifts as the solution of the NMP. This solution considers 3 shifts (for 3 technicians) in region 1, 1 shift in region 2, 1 shift in region 3 that includes a travel from site 6 to site 5, 1 shift in region 4 and 2 shifts in region 5.  $ \duration $ is a given set of possible shift durations, e.g. $ \{6h, 8h\} $. One reason for having multiple durations is related to the payment policy. Technicians should be paid for a minimum number of hours per shift. For example, if a technician works for any time less than 6 hours, they will be paid for the full 6 hours, while another technician working longer than 6 hours will be paid for 8 hours. In addition, access to the sites and the time spent at the sites should be within the SLA. Having multiple shift durations helps to avoid requesting unnecessary long access periods. 
    
    Migration costs include payments to the technicians and engineers. The network migration problem is to determine the order of upgrading the circuits in order to minimize the migration costs, i.e., building a set of minimum-cost plans for the technicians. 
    Table \ref{tab:param} presents the parameters and notations used in our model. 
    
    \begin{remark}
    Vectors, matrices, and scalars are represented by bold ($\boldsymbol{a}$), capital ($A$) and lower-case with regular font ($a$), respectively. We use $(.)^\top$ for the transpose operator, $|.|$ for the size of a set, $\boldsymbol{1}$ and $\boldsymbol{0}$ respectively, for vectors of 1s and 0s, adjusted to the required size. We use relational operators for element-wise comparison of two vectors.
    \end{remark}
    \begin{table}[t]
    \caption{Sets and parameters of the NMP}
    \centering
    \begin{tabular}{lp{13.5cm}}
        \toprule
        \textbf{\small Notation} & \textbf{\small Description}\\
        \midrule
        \textit{Sets}:\\
    	$ \siteSet $ & Set of sites (indexed by $ s $)\\
    	$ \sitePair $ & Set of site pairs $ \{s,s'\} $ with at least one circuit between them \\
    	$ \regSet $ & Set of regions (indexed by $ r $)\\
    	$ \siteSet_r $ & Set of sites in region $ r $\\
    	$\wSet$ & Set of available maintenance windows (indexed by $w$)\\
    	$ \cirSet $ & Set of circuits (indexed by $ c $)\\
    	$ \cirSet_{ss'} $ & Set of circuits between sites $ s $ and $ s'$\\
    	$ \cirSet_{r} $ & Set of circuits with at least one endpoint in the sites $ s\in \siteSet_{r} $\\
    	$ \duration $ & Ordered set of possible durations for a shift (e.g., 360 or 480 minutes)  with $\max(\duration)=\mwDur$, indexed by $ \delta $ \\[4mm]
    	\midrule
    	\textit{Parameters}:\\
    	$ \cirSS $ & Number of circuits between sites $ s $ and $ s' $ ($\cirSS = \phi_{s's} $)\\
    	$ \maxTechR $ & Max. number of technicians available in region $ r $ in a maintenance window\\
    	$ \maxCir $ & Max. number of circuits allowed to be migrated in a maintenance window\\
    	$ \maxEng $ & Max. number of engineers available in a maintenance window\\
    	$ \costT $ & Hourly cost of a technician\\
    	$ \costE $ & Hourly cost of an engineer\\
    	$\traTime_{ss'}$ & Travel time between sites $s$ and $s'$\\
        $\durmigr$ & Time required to migrate one endpoint of a circuit\\
        $\aEng$ & The number of technicians supported by one engineer\\[2mm]
        \bottomrule
    \end{tabular}
    \label{tab:param}
    \end{table}

\subsection{Problem Formulation}
\label{sec:opt_models}
We model the NMP as an ILP that returns a planning solution consisting of a set of shifts. The proposed formulation is amenable to the LBBD framework, meaning that we can decouple the problem into smaller subproblems that are easier to solve.

% ---------------------------------------------------------
% ---------------------------------------------------------

To begin with, assume that we have a set $\G $ of all possible shifts (in Section \ref{sec:CG} we \emph{implicitly} enumerate this set). For a shift $ \g \in \G $, the decision variable $ z_\g \in \mathbb{Z}_+ $ determines the number of times $ \g $ is assigned to the technicians. Every shift $ \g $ is characterized by ($i$) $ \durshift^{\g} $ its duration, ($ii$) $n_{ss'}^\g $ the number of circuit endpoints migrated between the pair of sites $\{s,s'\}$, and ($iii$) $ n_{\textsc{cir}}^\g $ the total number of migrated circuit endpoints in the shift. The total set of shifts is denoted by $ \G = \bigcup_{w\in \wSet}\G_w = \bigcup_{r\in \regSet, w\in\wSet}\G_{rw} $, where $ \G_w $ is the set of shifts  for a maintenance window $ w $ and $ \G_{rw} $ is the set of shifts for a technician located in region $ r $ during maintenance window $ w $. Denote by $m_{ss'w}\in \mathbb{Z}_+$,  a decision variable that determines the number of circuits between $\{s,s'\}$ migrated during $w$. The NMP is formulated as:
\begin{subequations}
\label{eq:NMP}
\begin{align}
\label{eq:obj_fun_model_II}
\min\ \ & \cost \sum_{\g \in \G}\durshift^{\g} z_{\g}\\
\text{s.t.}\ \ & \sum_{w\in \wSet} m_{ss'w} \geq \cirSS && \{s, s'\} \in \sitePair, s < s'
\label{eq: all_EPs_per_site_master_II} \\
&m_{ss'w} = m_{s'sw} && \{s, s'\} \in \sitePair, s < s', w \in \wSet
\label{eq: two_sites_same_EPs_per_MW_master_II} \\
&\sum_{\{s, s'\} \in \sitePair} m_{ss'w}  \leq 2\maxCir && w \in \wSet \label{eq: circuits_master_II}\\
& \sum_{\g \in \G_{w}}n^\g_{ss'} z_\g = m_{ss'w} && \{s, s'\} \in \sitePair, w \in \wSet  \label{eq: m_shift_eq_master_II} \\
& \sum_{\g \in \G_{rw}} z_\g \leq \maxTechR  && r\in \regSet, w \in \wSet  \label{eq: techs_in_r_mw_master_II} \\
&  \sum_{\g \in \G_{w}} z_\g \leq \aEng \maxEng  && w \in \wSet
    \label{eq: engineers_master_II} \\
& \boldsymbol{z} \in \mathbb{Z}_+^{|\G|},\ \boldsymbol{m} \in \mathbb{Z}_+^{|\sitePair|\times |\wSet|}, &&  \label{eq: domain_var_master_II}
\end{align}
\end{subequations}
where $\displaystyle\cost = \costT+\frac{\costE}{\aEng}$.
The objective function \eqref{eq:obj_fun_model_II} is the cost of the NMP, which is defined as the total technician and engineer costs over the duration of the migration. Constraints \eqref{eq: all_EPs_per_site_master_II} assure that all circuits between every two sites $s$ and $s'$ are migrated. Constraints \eqref{eq: two_sites_same_EPs_per_MW_master_II} enforce the number of migrated circuits from $s$ to $s'$ in $ w $ to be equal to the number of circuits migrated from $ s' $ to $ s $. Constraints \eqref{eq: circuits_master_II} establish the bound on the number of migrated circuits at every maintenance window. Through constraints  \eqref{eq: m_shift_eq_master_II},  variables $m_{ss'w}$ are determined by aggregating over the number of migrated circuits between $\{s,s'\}$ during the shifts at $w$. Constraints \eqref{eq: techs_in_r_mw_master_II} and \eqref{eq: engineers_master_II} ensure that at  $ w $, the number of available technicians and engineers are respected.  Constraints \eqref{eq: domain_var_master_II} define the variable domains.

The size of $\G$, the set of all shifts,  is an exponential function of the number of circuits, hence it is not reasonable (or even possible) to include them all in solving the model \eqref{eq:NMP}. Column generation is a method for implicitly enumerating such a large set of columns that relies on the duality theory for linear programming (LP).
In the presence of integer decision variables, the \emph{branch-and-price} (B\&P) algorithm combines the branch-and-bound framework for solving a MIP with the CG procedure. The performance of a B\&P  depends on the strength of the LP bound, as well as the employed branching and search strategy. It has been observed that branching on the variables of the master problem associated with the generated columns is not efficient and results in an unbalanced tree \citep{vance1998branch,vanderbeck2011branching}. 
Branching on the aggregate variables of the original formulations, in case of identical subproblems, is not typically sufficient to eliminate all fractional solutions. Although this branching scheme theoretically does not guarantee the integrality of the solution, it experimentally returns the integral solution for some instances. \cite{vanderbeck2011branching} proposes a generic branching scheme based on the aggregated value of the original variables when returning to non-identical systems. In our preliminary experiments though, specially for small to medium-size instances, the quality of the lower bound was  poor and did not improve adequately as the B\&P proceeded. Therefore, since branching on neither the subproblem nor the master problem variables for the NMP solved any of our instances, we concluded that a pure B\&P is not suitable to obtain exact solutions for the NMP.

In the next section, we develop an LBBD framework that entails multiple levels of decomposition, with a proof of optimality if one exists, or infeasibility when it does not.

\section{Solution Method} \label{sec:solution-method}
In what follows, we first discuss the general LBBD framework for the NMP, then we provide detailed discussions on each step of the algorithm.
\subsection{Logic-based Benders Decomposition}
In the problem formulation \eqref{eq:NMP}, decision variables $\boldsymbol{z}$ and $\boldsymbol{m}$ are only linked through \eqref{eq: m_shift_eq_master_II}. Furthermore,  except for the constraints \eqref{eq: all_EPs_per_site_master_II}, the problem is decomposable by maintenance windows. Using these two points, we propose the following decomposition. Denote by $\text{SP}^{\text{LBBD}}_w(\boldsymbol{m}_w)$, the problem of generating a  set of shifts with minimum cost for a given  $\boldsymbol{m}_w$. Let  $\eta_w$ be the  migration cost at  $w$.  We can reformulate the model \eqref{eq:NMP} as below:
\begin{equation}
\label{eq:NMP_BD_MP}
\min \ \Bigg\{ \sum_{w\in \wSet} \eta_w\ :\
 \eqref{eq: all_EPs_per_site_master_II} - \eqref{eq: circuits_master_II},\ \eta_w \geq \text{SP}^{\text{LBBD}}_w(\boldsymbol{m}_w),\ w \in \wSet,\ \boldsymbol{m} \in \mathbb{Z}_+^{|\sitePair|\times |\wSet|}\Bigg\}.
\end{equation}
% Note that, the subproblems $\text{SP}^{\text{LBBD}}_w(\boldsymbol{m}_w)$ are defined for each maintenance window, reducing the size of the problem we need to solve.
Formulation \eqref{eq:NMP_BD_MP} has the structure of a two-stage problem: the first-stage (master) problem decides on the number of migrated circuits between the site pairs at each maintenance window, along with an estimation on the cost of such a plan; the second-stage (recourse) problems verify if  it is feasible to migrate the assigned number of circuits with the available resources, and if so what is the actual cost of this migration. 

The \emph{Benders decomposition}  \citep{benders1962partitioning} is a well-established solution method for two-stage linear programs with continuous second-stage variables. After decomposing the problem into a master problem and a subproblem, it iteratively approximates the optimal solution to the recourse problem via Benders \emph{feasibility} and \emph{optimality} cuts (see Section \ref{sec:BD} for more details), derived using the LP duality theory. For two-stage problems with mixed-integer recourse decisions, logic-based
Benders decomposition \citep{hooker2003logic} is a generalization of the Benders decomposition (including its special case in the context of stochastic programming, the integer L-shaped method; see \cite{laporte1993integer,angulo2016improving}). For the NMP, the LBBD master problem is:
\begin{subequations}
\label{eq:NMP_LBBD_MP}
\begin{align}
\label{eq:obj_NMP_LBBD_MP}
\LBBDMP = \min \ \  & \sum_{w\in \wSet} \eta_w\\
\text{s.t.} \ \  & \eqref{eq: all_EPs_per_site_master_II} - \eqref{eq: circuits_master_II}\\
& (\eta_w,\boldsymbol{m}_w) \in \feascuts , && w \in \wSet \label{eq:-feasibilitycuts}\\
& (\eta_w,\boldsymbol{m}_w) \in \optcuts , && w \in \wSet \label{eq:optimality-cuts}\\
& \boldsymbol{m} \in \mathbb{Z}_+^{|\sitePair|\times |\wSet|}, \boldsymbol{\eta}\geq\boldsymbol{0}. && \label{eq: domain_NMP_LBBD_MP}
\end{align}
\end{subequations}
where $ \feascuts$ and $\optcuts $ are sets of feasibility and optimality cuts, respectively, and together represent an LBBD subproblem $\text{SP}^{\text{LBBD}}_w(\boldsymbol{m}_w)$. The LBBD starts by  $ \feascuts=\optcuts=\emptyset $ and at each iteration expands these sets with cuts if necessary until they are representative of the subproblem.
The LBBD relies on logical reasoning for obtaining feasibility and optimality cuts, and as long as we can have feasibility and optimality certificates, the LBBD subproblem can have any form. Having such a generic framework, the LBBD depends on the modeler for designing problem-specific cuts and, unlike the Benders decomposition, does not have a readily available cut development mechanism.

The  strength of the LBBD lies in the fact that it can integrate various optimization paradigms, most notably mixed-integer programming and constraint programming \citep{jain2001algorithms,hooker2012integrated}. 
Our first-stage problem $\LBBDMP$ is a variant of an assignment problem, which is suitable for a MIP-based solution approach. On the other hand, our second-stage problem takes $\overline{\boldsymbol{m}}_w$ as an argument and looks for the best plan for such an assignment. CP is an optimization paradigm particularly powerful for planning and scheduling problems.  
Accordingly, in our LBBD solution framework for the NMP, we have a MIP model as the master problem, and $|\wSet|$ many CP models as the subproblems.
In the remaining of this section, we first discuss the $\LBBDMP$ and develop two other levels of decomposition to make it more informed. Then, for a master problem solution $\overline{\boldsymbol{m}}_w$, we formulate the $\LBBDSP$ as a CP model, followed by the description of valid LBBD feasibility and optimality cuts that guarantee the convergence of the algorithm to an optimal solution.

% ---------------------------------------------------------
% ---------------------------------------------------------
% ---------------------------------------------------------
\subsection{The LBBD Master Problem: A Benders Decomposition} \label{sec:BD}
In its current form, $\LBBDMP$ is oblivious to the structure of the NMP. In this section, our goal is to make the solutions of $\LBBDMP$ more intelligent before passing them on to the $\LBBDSP$. In the literature of LBBD, it has been observed that adding a relaxation of the subproblem considerably improves the performance of the algorithm \citep{elci2020stochastic}. Often, this subproblem relaxation is in the form of an analytical expression based on the structure of the problem and is added as a bound to the master problem. We, however, resort to the LP relaxation of the subproblems to obtain valid inequalities for the $\LBBDMP$. 

By relaxing the integrality constraints of the subproblems, we now have integer first-stage and continuous second-stage decision variables and the new problem is amenable to the Benders decomposition, with $\BDMP = \LBBDMP$ as its master problem and the following subproblem for a maintenance window $w$:
\begin{subequations}
\label{eq:NMP_BD_SP}
\begin{align}
\label{eq:obj_NMP_BD_SP}
\BDSP = \min\ \  & \cost \sum_{\g \in \G_w}\durshift^{\g} z_{\g}\\
\text{s.t.}\ \  & \sum_{\g \in \G_{w}}n^\g_{ss'} z_\g \geq \overline{m}_{ss'w} && \{s, s'\} \in \sitePair, w \in \wSet  \label{eq: m_shift_eq_NMP_BD_SP} \\
& \eqref{eq: techs_in_r_mw_master_II} - \eqref{eq: engineers_master_II} \\
& \boldsymbol{z} \geq\boldsymbol{0}. && \label{eq: domain_NMP_BD_SP}
\end{align}
\end{subequations}
Considering that \eqref{eq:obj_NMP_BD_SP} has an exponential number of variables, in the next section we apply Dantzig-Wolfe decomposition principles to develop a column generation procedure for systematically adding them to the set of columns. For a review on the Dantzig-Wolfe decomposition and column generation, unfamiliar reader may refer to  \cite{chvatal1983linear}.

\subsubsection{The Benders Subproblem: A Column Generation Method.}\label{sec:CG}
A CG solution method starts with solving the \emph{restricted master problem}, defined as the original problem with a (potentially empty) subset of all the columns. Then, the optimal dual solutions are passed to a pricing problem (the CG subproblem) that checks their feasibility in the LP dual of the original problem, and if not, adds an improving column to the master problem. Once the pricing problem determines that a feasible dual solution is found, the CG stops as we have reached the optimality. 
By design, the pricing problem implicitly considers all the columns by using the properties that define valid columns of the master problem. Although the Dantzig-Wolfe decomposition lays out a precise scheme for decomposing a problem for the CG method, the master and pricing problems are often built by problem-specific modeling practices. In defining the ``shifts'' and $z_\g$ as the decision variables we have used such an approach which can directly be translated into the pricing problems.  

\textit{The CG master problem.} Master problem of the NMP selects the best set of shifts for the technicians among a subset of columns $\G'\subseteq\G$, and its pricing problems generate improving shifts. The CG master problem is:
\begin{subequations}
\label{eq:NMP_CG_RMP}
\begin{align}
\label{eq:obj_NMP_CG_RMP}
\CGMP = \min\ \  & \cost \sum_{\g \in \Gp}\durshift^{\g} z_{\g} \\
\text{s.t.}\ \ & \sum_{\g \in \Gp_{w}}n^\g_{ss'} z_\g \geq \overline{m}_{ss'w} && \{s, s'\} \in \sitePair  \label{eq: m_shift_eq_NMP_CG_RMP} \\
& \sum_{\g \in \Gp_{rw}} z_\g \leq \maxTechR  && r\in \regSet  \label{eq: techs_in_r_mw_NMP_CG_RMP} \\
&  \sum_{\g \in \Gp_{w}} z_\g \leq \aEng \maxEng  && 
    \label{eq: engineers_NMP_CG_RMP} \\
& \boldsymbol{z} \geq \boldsymbol{0}.  && \label{eq: domain_NMP_CG_RMP}
\end{align}
\end{subequations}
Constraints of $\CGMP$ correspond to the constraints \eqref{eq: m_shift_eq_master_II} - \eqref{eq: domain_var_master_II}, except for the set of columns $\G'$, and the fact that they are for a single maintenance window $w$.

\textit{The CG pricing problem.} In the definition of the CG subproblems, for the sake of brevity, we drop the index $ \g $ from the decision variables. 
Denote by $\overline{\boldsymbol{\dual}}^{\eqref{eq: m_shift_eq_NMP_CG_RMP}}, \overline{\boldsymbol{\dual}}^{\eqref{eq: techs_in_r_mw_NMP_CG_RMP}},  \overline{\boldsymbol{\dual}}^{\eqref{eq: engineers_NMP_CG_RMP}}$, the optimal dual solutions associated with constraints ${\eqref{eq: m_shift_eq_NMP_CG_RMP}}$, ${\eqref{eq: techs_in_r_mw_NMP_CG_RMP}}$ and ${\eqref{eq: engineers_NMP_CG_RMP}} $, respectively. The pricing problem generating a shift for region $ r $ and maintenance window $ w $ is as follows:
\begin{subequations}
\label{eq:pp}
\begin{align}
\CGSP = \min \ \  &\cost\durshift -\sum_{s\in \siteSet_r}\sum_{s'\in \siteSet} n_{ss'}\overline{\dual}_{ss'}^{\eqref{eq: m_shift_eq_NMP_CG_RMP}} - \overline{\dual}^{\eqref{eq: techs_in_r_mw_NMP_CG_RMP}}_{rw} - \overline{\dual}^{\eqref{eq: engineers_NMP_CG_RMP}}_{w} \label{eq:pp-obj}\\
\text{s.t.}\ \  & \{\text{Constraints defining a valid shift}\}\label{eq:validshift}\\
& n_{ss'}\in \mathbb{Z}_+ \qquad&& s \in \siteSet_r, \{s,s'\}\in\sitePair\\
& \durshift \geq 0.
\end{align}
\end{subequations}
Objective function \eqref{eq:pp-obj} is the reduced-cost that determines if $\CGMP$ is at optimality. If not, the pricing problem generates a column that corresponds to a ``valid'' shift, and is characterized by its duration $\durshift$ and the number of circuits migrated between each site pair. A shift is made of a sequence of site visits by the technician. So constraints \eqref{eq:validshift} define a valid shift as a connected path over the sites in the region $r$ such that the duration of the shift does not exceed the maximum possible duration of the maintenance window. Furthermore, as the endpoints of a circuit should be migrated by two technicians in the same maintenance window, at most one of the endpoints of every circuit $ c\in \cirSet_{ss'} $ can be migrated in a  shift.  In Section \ref{sec:CP-model}, we present a constraint programming model for generating a set of valid shifts in a maintenance window $w$, for all the regions and their technicians. Our preliminary experiments revealed that, because of the routing decisions in the pricing problems and the presence of loop elimination constraints, $\CGSP$ as a CP model performs much better than a MIP. Therefore, in lieu of model \eqref{eq:pp}, we fix the region in the CP model of Section \ref{sec:CP-model} and solve it for one technician, with the objective function \eqref{eq:pp-obj}. 
\begin{remark}\label{rem:hybriddCG}
To accelerate the solution process, we use a hybrid CG, where first an auxiliary pricing problem  generates improving columns with ``ordered paths'' that only have $ (s,s'), s < s'$ links (model description is given in the e-companion). After the auxiliary problem converges, we solve $\CGSP$ to verify the optimality. We observed that this two-subproblem strategy greatly improves the performance of the algorithm. The reason is that on backbone networks, it is quite possible that many regions have a few number of sites. If a region has up to two sites, the auxiliary problem alone guarantees the optimality.
\end{remark}

Let $\CGSPSol$ be the optimal solution of the $\CGSP$. If $\CGSPSol\geq 0, \forall r\in\regSet$, then the CG procedure stops. Otherwise, for each $r\in\regSet$ with $\CGSPSol< 0$ we add the generated column to $\G'$ and repeat the process.
As a result of having the constraints \eqref{eq: m_shift_eq_NMP_CG_RMP}, we require an initial set of columns $\G'\neq\emptyset$ that make the $\CGMP$ feasible. The usual approach for generating such $\G'$ is to go through an initial phase (\texttt{INIT}) where an artificial non-negative decision variable $\artvar$ is added to each ``$\geq$'' constraint with a positive right-hand-side, and the objective function is replaced with $\boldsymbol{1}^\top\boldsymbol{\rho}$. If the CG procedure for the new problem stops with an optimal value equal to zero, then the artificial decision variables are removed from the problem, the original objective function is brought back and the generated columns are selected as $\G'$. Otherwise, if the optimal objective value is positive, we can conclude that the original problem is infeasible.
% ---------------------------------------------------------
% ---------------------------------------------------------
% ---------------------------------------------------------
\subsubsection{Benders Cuts.} Depending on the status of the $\CGMP$ after solving, we might need to add feasibility (optimality) cuts to $\feascuts$ ($\optcuts$). Benders decomposition provides us with off-the-shelf cuts through the dual solutions of the subproblems. For our problem, because the subproblems are solved via CG,  the feasibility cuts are not immediately clear. Next, we present the Benders feasibility and optimality cuts for the NMP and show that they are valid, despite being obtained from a restricted set of columns in the $\CGSP$.

\textit{Benders feasibility cut}.
Assume that, at the end of the \texttt{INIT} phase,  the CG procedure stops  with $\boldsymbol{\artvar} \neq \boldsymbol{0}$, and $\overline{\boldsymbol{\dual}}^{\eqref{eq: m_shift_eq_NMP_CG_RMP}}, \overline{\boldsymbol{\dual}}^{\eqref{eq: techs_in_r_mw_NMP_CG_RMP}},  \overline{\dual}^{\eqref{eq: engineers_NMP_CG_RMP}}$ are returned from the (modified) $\CGMP$. The Benders feasibility cut to be added to $\feascuts$ is as follows:
\begin{equation}\label{eq:BD_feas_cut}
    0 \geq \sum_{\{s, s'\} \in \sitePair} \overline{\dual}^{\eqref{eq: m_shift_eq_NMP_CG_RMP}}_{ss'} m_{ss'w} + \sum_{r \in \regSet} \overline{\dual}^{\eqref{eq: techs_in_r_mw_NMP_CG_RMP}}_{r}\maxTechR + \overline{\dual}^{\eqref{eq: engineers_NMP_CG_RMP}}\aEng \maxEng
\end{equation}
To show that \eqref{eq:BD_feas_cut} is indeed a feasibility cut for the $\BDMP$ and cuts off the infeasible solution $\overline{\boldsymbol{m}}_w$, we first prove that $(\overline{\boldsymbol{\dual}}^{\eqref{eq: m_shift_eq_NMP_CG_RMP}}, \overline{\boldsymbol{\dual}}^{\eqref{eq: techs_in_r_mw_NMP_CG_RMP}},  \overline{\dual}^{\eqref{eq: engineers_NMP_CG_RMP}})$ constitutes a certificate of infeasibility for the $\BDMP$, even if it is derived from the $\CGMP$ with $\Gp\subseteq\G$. In the following theorem, we show that this  is true for any Dantizig-Wolfe decomposition at the end of the \texttt{INIT} phase.
\begin{theorem}\label{eq:thm_BD_feas_cut}
Consider $\overline{P}(\mathcal{I}')$, the restricted master problem of a Dantzig-Wolfe decomposition at the end of the \texttt{INIT} phase, and $P(\mathcal{I}',c_\mathcal{I}')$, the problem that is obtained by removing the artificial variables $\boldsymbol{\artvar}$ from $\overline{P}(\mathcal{I}')$ and bringing back the original objective function:

\begin{minipage}{.49\textwidth}
\raggedright
    \begin{align*}
        \overline{P}(\mathcal{I}') = \min\ \  & {\boldsymbol{1}^\top} \boldsymbol{\artvar}\\
        \text{s.t.}\ \  & A_{\mathcal{I'}}\boldsymbol{x}_{\mathcal{I}'} + \boldsymbol{\artvar} \geq \boldsymbol{b} \\
        & \boldsymbol{x}_{\mathcal{I}'}, \boldsymbol{\artvar} \geq \boldsymbol{0},\\
    \end{align*}
    
\end{minipage}% <---------------- Note the use of "%"
\begin{minipage}{.49\textwidth}
\raggedleft
    \begin{align*}
        P(\mathcal{I}',\boldsymbol{c}_{\mathcal{I}'}) = \min\ \  & {\boldsymbol{c}_{\mathcal{I}'}^\top} \boldsymbol{x}_{\mathcal{I}'}\\
        \text{s.t.}\ \  & A_{\mathcal{I'}}\boldsymbol{x}_{\mathcal{I}'}\geq \boldsymbol{b} \\
        & \boldsymbol{x}_{\mathcal{I}'} \geq \boldsymbol{0},\\
    \end{align*}
    
\end{minipage}
where $ \emptyset \neq  \mathcal{I}'\subseteq\mathcal{I}$,  and $\boldsymbol{c}_{\mathcal{I}'}, A_{\mathcal{I}'}$ are the cost vector and columns associated with $\mathcal{I}'$. 
$\boldsymbol{\pi}$ is a certificate of infeasibility for $P(\mathcal{I}',\boldsymbol{c}_{\mathcal{I}'})$ only if it is a certificate of infeasibility for $P(\mathcal{I},\boldsymbol{c}_{\mathcal{I}})$. 
\end{theorem}
\proofNoH{
We first find the certificate of infeasibility for $P(\mathcal{I}',\boldsymbol{c}_{\mathcal{I}'})$. Then we show that it is also a certificate of infeasibility for $P(\mathcal{I},\boldsymbol{c}_{\mathcal{I}})$. Let $(\overline{\boldsymbol{x}}_{\mathcal{I}'},\overline{\boldsymbol{\artvar}})$ denote the optimal solution of $\overline{P}(\mathcal{I}')$.
\begin{enumerate}[leftmargin=*]
    \item From (a variant of) the Farkas Lemma we know that exactly one of the following system of inequalities has a solution \citep{matousek2007understanding}:
$$\text{(I) } A_{\mathcal{I'}}\boldsymbol{x}_{\mathcal{I}'}\geq \boldsymbol{b}, \boldsymbol{x}\geq\boldsymbol{0}, \qquad\quad \text{(II) } \boldsymbol{\pi}^\top A_{\mathcal{I'}}\leq \boldsymbol{0}, \boldsymbol{\pi}^\top\boldsymbol{b} > 0, \boldsymbol{\pi}\geq \boldsymbol{0}. $$
If $\overline{\boldsymbol{\artvar}}\neq\boldsymbol{0}$, then $\{\boldsymbol{x} :  A_{\mathcal{I'}}\boldsymbol{x}\geq \boldsymbol{b}, \boldsymbol{x}\geq\boldsymbol{0}\} = \emptyset$ and $P(\mathcal{I}',\boldsymbol{c}_{\mathcal{I}'})$ is infeasible. Therefore, there exits a $\boldsymbol{\pi}_0$ that satisfies the inequalities of (II) and is a certificate of infeasibility for $P(\mathcal{I}',\boldsymbol{c}_{\mathcal{I}'})$. Consider the dual of $\overline{P}(\mathcal{I}')$ as follows:
$$\overline{D}(\mathcal{I}') = \max \{{\boldsymbol{\pi}^\top} \boldsymbol{b}: \boldsymbol{\pi}^\top A_{\mathcal{I'}} \leq \boldsymbol{0}^\top, \boldsymbol{\pi} \leq \boldsymbol{1}, \boldsymbol{\pi} \geq \boldsymbol{0}\}.$$
From $\overline{\boldsymbol{\artvar}}\neq 0$, we have $\boldsymbol{1}^\top \overline{\boldsymbol{\artvar}} = \overline{\boldsymbol{\pi}}^\top \boldsymbol{b} > 0$, with $\overline{\boldsymbol{\pi}}$ the optimal solution of $\overline{D}(\mathcal{I}')$. Clearly $\overline{\boldsymbol{\pi}}$ satisfies (II) and we can set $\boldsymbol{\pi}_0 = \overline{\boldsymbol{\pi}}$. So $\overline{\boldsymbol{\pi}}$ is the Farkas certificate of $P(\mathcal{I}',\boldsymbol{c}_{\mathcal{I}'})$.
\item Let $D(\mathcal{I},\boldsymbol{0}) = \max\{\boldsymbol{\pi}^\top\boldsymbol{b}: \boldsymbol{\pi}^\top A_{\mathcal{I}}\leq \boldsymbol{0}^\top, \boldsymbol{\pi}\geq \boldsymbol{0} \}$ be the dual of $P(\mathcal{I},\boldsymbol{0})$. Since $\overline{P}(\mathcal{I}')$ is the last restricted master problem at phase one, no improving column with negative reduced-cost is found through the pricing problem with the dual solution of $\overline{P}(\mathcal{I}')$, i.e., $\boldsymbol{0}^\top - \overline{\boldsymbol{\pi}}^\top A_{\mathcal{I}} \geq 0$ meaning that $\overline{\boldsymbol{\pi}}$ is feasible for $D(\mathcal{I},\boldsymbol{0})$. 

\item As $\overline{\boldsymbol{\artvar}}\neq\boldsymbol{0}$, we have $\{\boldsymbol{x} :  A_{\mathcal{I}}\boldsymbol{x}\geq \boldsymbol{b}, \boldsymbol{x}\geq\boldsymbol{0}\} = \emptyset$, $P(\mathcal{I},\boldsymbol{c}_{\mathcal{I}})$ is infeasible, and so is $P(\mathcal{I},\boldsymbol{0})$. Accordingly, using the Farkas Lemma one more time, there exists a $\hat{\boldsymbol{\pi}}$ such that $\hat{\boldsymbol{\pi}}^\top A_{\mathcal{I}}\leq \boldsymbol{0}, \hat{\boldsymbol{\pi}}^\top\boldsymbol{b} > 0, \hat{\boldsymbol{\pi}}\geq \boldsymbol{0} $. Consider the ray $\overline{\boldsymbol{\pi}} + \lambda \hat{\boldsymbol{\pi}}, \lambda\geq 0$. Then:
\begin{eqnarray*}
    (\overline{\boldsymbol{\pi}} + \lambda \hat{\boldsymbol{\pi}})^\top A_{\mathcal{I}} = \overbrace{\overline{\boldsymbol{\pi}}^\top A_{\mathcal{I}}}^{
    \leq \boldsymbol{0}} + \lambda \overbrace{\hat{\boldsymbol{\pi}}^\top A_{\mathcal{I}}}^{\leq\boldsymbol{0}} &\Longrightarrow & (\overline{\boldsymbol{\pi}} + \lambda \hat{\boldsymbol{\pi}})^\top A_{\mathcal{I}} \leq \boldsymbol{0}, \\
    (\overline{\boldsymbol{\pi}} + \lambda \hat{\boldsymbol{\pi}})^\top \boldsymbol{b} = \overbrace{\overline{\boldsymbol{\pi}}^\top \boldsymbol{b}}^{
    > \boldsymbol{0}} + \lambda \overbrace{\hat{\boldsymbol{\pi}}^\top \boldsymbol{b}}^{>\boldsymbol{0}} &\Longrightarrow & (\overline{\boldsymbol{\pi}} + \lambda \hat{\boldsymbol{\pi}})^\top \boldsymbol{b} > \boldsymbol{0}.
\end{eqnarray*}
Therefore, $\overline{\boldsymbol{\pi}} + \lambda \hat{\boldsymbol{\pi}}$ is feasible for $D(\mathcal{I},\boldsymbol{0})$ with a positive objective value $(\overline{\boldsymbol{\pi}} + \lambda \hat{\boldsymbol{\pi}})^\top \boldsymbol{b}$. As $\lambda \rightarrow +\infty$, so does $(\overline{\boldsymbol{\pi}} + \lambda \hat{\boldsymbol{\pi}})^\top \boldsymbol{b}$, proving that  $\overline{\boldsymbol{\pi}}$ is a certificate of unboundedness for $D(\mathcal{I},\boldsymbol{0})$, hence a certificate of infeasibility for $P(\mathcal{I},\boldsymbol{0})$ and $P(\mathcal{I},\boldsymbol{c}_{\mathcal{I}})$.\hfill\Halmos\\[-2mm]
\end{enumerate}}

\noindent From the  above discussions, $(\overline{\boldsymbol{\dual}}^{\eqref{eq: m_shift_eq_NMP_CG_RMP}}, \overline{\boldsymbol{\dual}}^{\eqref{eq: techs_in_r_mw_NMP_CG_RMP}},  \overline{\dual}^{\eqref{eq: engineers_NMP_CG_RMP}})$ is a proof of infeasibility for $\BDMP$ with  $0 <\sum_{\{s, s'\} \in \sitePair} \overline{\pi}^{\eqref{eq: m_shift_eq_NMP_CG_RMP}}_{ss'} \overline{m}_{ss'w} + \sum_{r \in \regSet} \overline{\dual}^{\eqref{eq: techs_in_r_mw_NMP_CG_RMP}}_{w}\maxTechR + \overline{\dual}^{\eqref{eq: engineers_NMP_CG_RMP}}\aEng\maxEng$ and can be removed using the inequality  \eqref{eq:BD_feas_cut}. Here we should mention that, we can also add the LBBD feasibility cuts (Section \ref{sec:feas-lbbd-cut}) for cutting infeasible solutions. However, as will be discussed later, our LBBD feasibility cuts are costly and  inequalities \eqref{eq:BD_feas_cut} improve the overall efficiency of the method.

\textit{Benders optimality cut}. If at the end of the \texttt{INIT} phase, $\boldsymbol{\rho} = \boldsymbol{0}$, then $\CGSP$ is feasible. Now  we should check if $\overline{\eta}_w$ is an accurate estimate of the $\BDSPSol$. Because the CG procedure stops at  optimality, we can treat the missing columns in $\G'$ as non-basic variables. Therefore it is clear that the optimal dual solutions from solving $\CGMP$ with $\G'$ are the same as the optimal dual solutions of $\BDSP$.   If  $\overline{\eta}_w \geq \CGSPSol$, then $(\overline{\eta}_w,\overline{\boldsymbol{m}}_w)$ is feasible in $\BDMP$. Otherwise, we cut it off by adding the following inequality to $\optcuts$:
\begin{equation}\label{eq:BD-opt-cut}
    {\eta}_w \geq \CGSPSol- \sum_{\{s, s'\} \in \sitePair} \overline{\dual}^{\eqref{eq: m_shift_eq_NMP_CG_RMP}}_{ss'} (\overline{m}_{ss'w} - m_{ss'w}).
\end{equation}

\subsection{The LBBD Subproblem: A Constraint Programming Model}\label{sec:CP-model}
With its roots in logic, CP is a modeling framework for combinatorial problems and has proven quite powerful for making planning, scheduling and routing decisions \citep{cire2016logic}. Integration of CP and mathematical programming models through LBBD \citep{hooker2000scheme}, and CG \citep{rousseau2004solving}, often results in solution algorithms that outperform methods relying solely on either of the two.  \cite{hooker2018constraint} provide a review on  integration of CP and Operations Research. 
In this section, considering the presence of planning and routing decisions, we formulate the $\LBBDSP$ as a CP model.

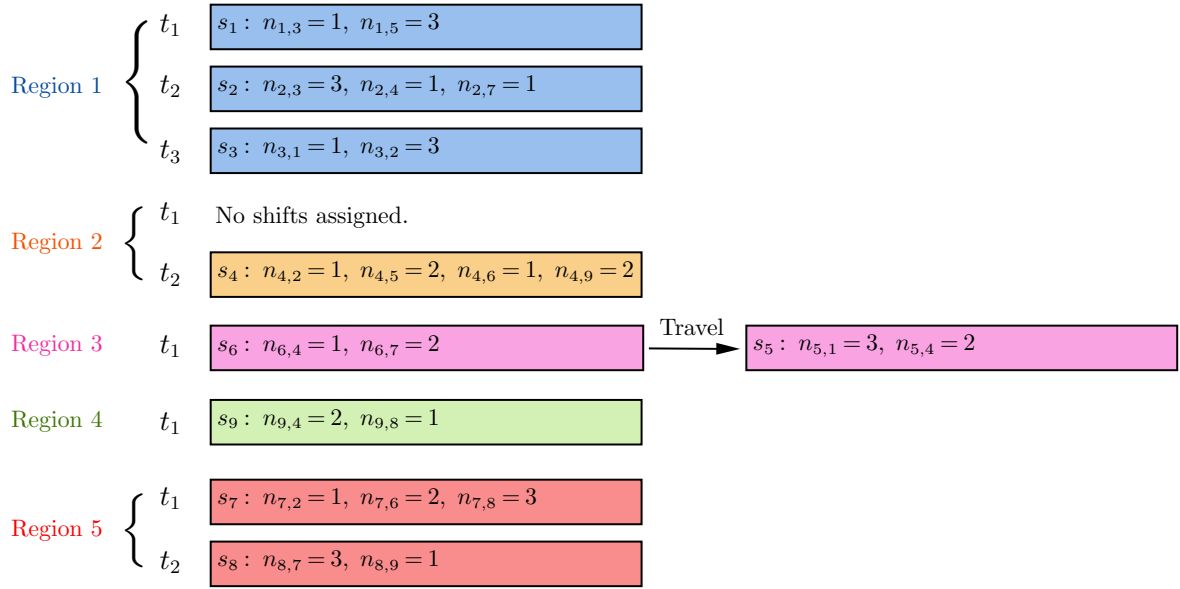
\begin{figure}[t]
    	\centering
    	\tikzset{every picture/.style={line width=0.75pt}} %set default line width to 0.75pt        

\begin{tikzpicture}[x=0.75pt,y=0.75pt,yscale=-1,xscale=1]
%uncomment if require: \path (0,372); %set diagram left start at 0, and has height of 372

%Shape: Rectangle [id:dp34638254425937665] 
\draw  [fill={rgb, 255:red, 74; green, 144; blue, 226 }  ,fill opacity=0.57 ] (108,22) -- (325.32,22) -- (325.32,44.5) -- (108,44.5) -- cycle ;
%Shape: Rectangle [id:dp787233430042596] 
\draw  [fill={rgb, 255:red, 74; green, 144; blue, 226 }  ,fill opacity=0.57 ] (108,53.25) -- (325.32,53.25) -- (325.32,75.75) -- (108,75.75) -- cycle ;
%Shape: Rectangle [id:dp5404671421604774] 
\draw  [fill={rgb, 255:red, 74; green, 144; blue, 226 }  ,fill opacity=0.57 ] (108,84.5) -- (325.32,84.5) -- (325.32,107) -- (108,107) -- cycle ;
%Shape: Rectangle [id:dp22157553438171218] 
\draw  [fill={rgb, 255:red, 241; green, 8; blue, 178 }  ,fill opacity=0.37 ] (378.15,184) -- (595.46,184) -- (595.46,206.5) -- (378.15,206.5) -- cycle ;
%Shape: Rectangle [id:dp21341415568780875] 
\draw  [fill={rgb, 255:red, 241; green, 8; blue, 178 }  ,fill opacity=0.37 ] (108,184) -- (326.15,184) -- (326.15,206.5) -- (108,206.5) -- cycle ;
%Shape: Rectangle [id:dp8396698816144343] 
\draw  [fill={rgb, 255:red, 245; green, 166; blue, 35 }  ,fill opacity=0.53 ] (108,147) -- (325.32,147) -- (325.32,169.5) -- (108,169.5) -- cycle ;
%Shape: Rectangle [id:dp670354314888821] 
\draw  [fill={rgb, 255:red, 184; green, 233; blue, 134 }  ,fill opacity=0.61 ] (108,221.5) -- (325.32,221.5) -- (325.32,244) -- (108,244) -- cycle ;
%Shape: Rectangle [id:dp8406435466572804] 
\draw  [fill={rgb, 255:red, 244; green, 33; blue, 33 }  ,fill opacity=0.51 ] (108,261.75) -- (325.32,261.75) -- (325.32,284.25) -- (108,284.25) -- cycle ;
%Shape: Rectangle [id:dp5216383930367231] 
\draw  [fill={rgb, 255:red, 244; green, 33; blue, 33 }  ,fill opacity=0.51 ] (108,293) -- (325.32,293) -- (325.32,315.5) -- (108,315.5) -- cycle ;
%Straight Lines [id:da86072766511914] 
\draw    (329,195.5) -- (365,195.74) ;
\draw [shift={(375,195.75)}, rotate = 180.28] [fill={rgb, 255:red, 0; green, 0; blue, 0 }  ][line width=0.08]  [draw opacity=0] (12,-3) -- (0,0) -- (12,3) -- cycle    ;

% Text Node
\draw (6,57) node [anchor=north west][inner sep=0.75pt]  [font=\small,color={rgb, 255:red, 7; green, 77; blue, 159 }  ,opacity=1 ,xscale=0.9,yscale=0.9] [align=left] {{\small Region 1}};
% Text Node
\draw (61,28.4) node [anchor=north west][inner sep=0.75pt]  [font=\Large,xscale=0.9,yscale=0.9]  {$\Bigg\{$};
% Text Node
\draw (6,136) node [anchor=north west][inner sep=0.75pt]  [font=\small,color={rgb, 255:red, 245; green, 80; blue, 0 }  ,opacity=1 ,xscale=0.9,yscale=0.9] [align=left] {{\small Region 2}};
% Text Node
\draw (6,187.25) node [anchor=north west][inner sep=0.75pt]  [font=\small,color={rgb, 255:red, 241; green, 45; blue, 160 }  ,opacity=1 ,xscale=0.9,yscale=0.9] [align=left] {{\small Region 3}};
% Text Node
\draw (6,225.5) node [anchor=north west][inner sep=0.75pt]  [font=\small,color={rgb, 255:red, 65; green, 117; blue, 5 }  ,opacity=1 ,xscale=0.9,yscale=0.9] [align=left] {{\small Region 4}};
% Text Node
\draw (6,280.5) node [anchor=north west][inner sep=0.75pt]  [font=\small,color={rgb, 255:red, 249; green, 0; blue, 0 }  ,opacity=1 ,xscale=0.9,yscale=0.9] [align=left] {{\small Region 5}};
% Text Node
\draw (61,122.4) node [anchor=north west][inner sep=0.75pt]  [font=\Large,xscale=0.9,yscale=0.9]  {$\Big\{$};
% Text Node
\draw (61,267.4) node [anchor=north west][inner sep=0.75pt]  [font=\Large,xscale=0.9,yscale=0.9]  {$\Big\{$};
% Text Node
\draw (81,25.4) node [anchor=north west][inner sep=0.75pt]  [xscale=0.9,yscale=0.9]  {$t_{1}$};
% Text Node
\draw (81,120.4) node [anchor=north west][inner sep=0.75pt]  [xscale=0.9,yscale=0.9]  {$t_{1}$};
% Text Node
\draw (81,187.4) node [anchor=north west][inner sep=0.75pt]  [xscale=0.9,yscale=0.9]  {$t_{1}$};
% Text Node
\draw (81,226.4) node [anchor=north west][inner sep=0.75pt]  [xscale=0.9,yscale=0.9]  {$t_{1}$};
% Text Node
\draw (81,264.4) node [anchor=north west][inner sep=0.75pt]  [xscale=0.9,yscale=0.9]  {$t_{1}$};
% Text Node
\draw (81,56.4) node [anchor=north west][inner sep=0.75pt]  [xscale=0.9,yscale=0.9]  {$t_{2}$};
% Text Node
\draw (81,151.4) node [anchor=north west][inner sep=0.75pt]  [xscale=0.9,yscale=0.9]  {$t_{2}$};
% Text Node
\draw (81,296.4) node [anchor=north west][inner sep=0.75pt]  [xscale=0.9,yscale=0.9]  {$t_{2}$};
% Text Node
\draw (81,89.4) node [anchor=north west][inner sep=0.75pt]  [xscale=0.9,yscale=0.9]  {$t_{3}$};
% Text Node
\draw (110,25.4) node [anchor=north west][inner sep=0.75pt]  [font=\small,xscale=0.9,yscale=0.9]  {$s_{1} :\ n_{1,3} =1,\ n_{1,5} =3$};
% Text Node
\draw (110,56.65) node [anchor=north west][inner sep=0.75pt]  [font=\small,xscale=0.9,yscale=0.9]  {$s_{2} :\ n_{2,3} =3,\ n_{2,4} =1,\ n_{2,7} =1$};
% Text Node
\draw (110,87.9) node [anchor=north west][inner sep=0.75pt]  [font=\small,xscale=0.9,yscale=0.9]  {$s_{3} :\ n_{3,1} =1,\ n_{3,2} =3$};
% Text Node
\draw (109,122.5) node [anchor=north west][inner sep=0.75pt]  [xscale=0.9,yscale=0.9] [align=left] {{\small No shifts assigned.}};
% Text Node
\draw (110,150.4) node [anchor=north west][inner sep=0.75pt]  [font=\small,xscale=0.9,yscale=0.9]  {$s_{4} :\ n_{4,2} =1,\ n_{4,5} =2,\ n_{4,6} =1,\ n_{4,9} =2$};
% Text Node
\draw (110,187.9) node [anchor=north west][inner sep=0.75pt]  [font=\small,xscale=0.9,yscale=0.9]  {$s_{6} :\ n_{6,4} =1,\ n_{6,7} =2$};
% Text Node
\draw (333,179) node [anchor=north west][inner sep=0.75pt]  [font=\small,xscale=0.9,yscale=0.9] [align=left] {{\small Travel}};
% Text Node
\draw (380.15,187.4) node [anchor=north west][inner sep=0.75pt]  [font=\small,xscale=0.9,yscale=0.9]  {$s_{5} :\ n_{5,1} =3,\ n_{5,4} =2$};
% Text Node
\draw (110,224.9) node [anchor=north west][inner sep=0.75pt]  [font=\small,xscale=0.9,yscale=0.9]  {$s_{9} :\ n_{9,4} =2,\ n_{9,8} =1$};
% Text Node
\draw (110,265.15) node [anchor=north west][inner sep=0.75pt]  [font=\small,xscale=0.9,yscale=0.9]  {$s_{7} :\ n_{7,2} =1,\ n_{7,6} =2,\ n_{7,8} =3$};
% Text Node
\draw (110,296.4) node [anchor=north west][inner sep=0.75pt]  [font=\small,xscale=0.9,yscale=0.9]  {$s_{8} :\ n_{8,7} =3,\ n_{8,9} =1$};

\end{tikzpicture}
    	\caption{A feasible plan for the solution presented in Figure \ref{fig:configuration}}
    	\label{fig:configuration-plan}
    \end{figure}

For a maintenance window $w$, a CP model serving as our LBBD subproblem creates a \emph{plan}, i.e., a set of shifts corresponding to a set of technicians working during the maintenance window $w$. For example, assume that the solution illustrated in Figure \ref{fig:configuration} is planned for one maintenance window. A plan for this solution consists of determining for each technician, a connected path along with the number of circuit endpoints migrated between each site pair. One possible such plan is depicted in Figure \ref{fig:configuration-plan}, where a set of boxes in front of a technician $t$ shows its shift for the current maintenance window. Every box includes a site $s$ and the number of circuit endpoints $n_{ss'}$ that $t$ migrates in this shift. We see that in the CP model we need to treat the technicians as individual entities, unlike the engineers and circuit endpoints. Therefore, for each region $r$ at maintenance window $w$, we define $\T_{rw} = \{1,\dots,\maxTechR\}$ as its set of available technicians, and the CP model generates shifts for the technicians that belong to $\bigcup_{r\in\regSet}\T_{rw}$. It is clear that the newly introduced technician symmetry stays within the CP model, since the only link between $\LBBDMP$ and its subproblems is through the $\boldsymbol{m}$ decision variables defined for site pairs and independent of the technician (see Section \ref{sec:LBBD_cuts}).

Next, we describe the variable types, functions and global constraints used in the CP model, followed by the model description. While we use names and conventions of the CP Optimizer \citep{laborie2018ibm}, equivalent notions exist for other notable CP solvers such as the CP-SAT from Google OR Tools \citep{ortools}.\\[-8mm]
\begin{longtable}{l p{13cm}}
    \hspace{-1mm}- $\IntegerVar()$ & An \emph{integer} decision variable. If defined as \emph{Optional}, it can be absent from the CP solution.\\ 
    \hspace{-1mm}- $\IntervalVar()$ & An \emph{interval} decision variable, modeling a time interval characterized by start, end and length. It can be optional.\\ 
    \hspace{-1mm}- $\SequenceVar(\mathcal{X})$ & A \emph{sequence} decision variable defining a total order over a set of interval variables $\mathcal{X}$. If all members of $\mathcal{X}$ are absent, the sequence becomes empty.\\
    \hspace{-1mm}- $\PresenceOf (a)$ & Returns 1 if $a$ is present in the solution, 0 otherwise.\\ 
    \hspace{-1mm}- $\EndOf(a)$ & Returns the finish time of the task if the interval variable $a$ is present in the solution.\\ 
    \hspace{-1mm}- $\LengthOf(a)$ & Returns the time spent on the task if the interval variable $a$ is present in the solution.\\
    \hspace{-1mm}- $\Span(a,\mathcal{X})$ & Enforces the interval variable $a$ to span the set of interval variables $\mathcal{X}$.\\
    \hspace{-1mm}- $\NoOverlap (\ell,\travel)$ & Enforces the intervals in the sequence $\ell$ to be disjoint, while taking into account the transition times $T$.\\
    \hspace{-1mm}- $\IfThen(p,q)$ & Implements the logical constraint $ p \rightarrow q$.\\[-5mm]
\end{longtable}
In addition to the sets and parameters described in Section \ref{sec:prob}, we need the variables defined in Table \ref{tab:CP_var}.
\begin{table}[t]
\centering
\caption{Decision variables of the constraint programming model}
\label{tab:CP_var}
\begin{tabular}{lllcp{6.4cm}}
\toprule
\textbf{\small Variable} & \textbf{\small Type} & \textbf{\small Initial domain} & \textbf{\small Optional} & \textbf{\small Description}\\
\midrule
    $\techsiteTS $ & $\IntervalVar$ & $ [[0,\mwDur],[0,\mwDur]] $ & $\checkmark$ &  Time plan for tech $t$ of $r$ at site $s$\\
    $\seqT$ & $\SequenceVar$ & $\mathfrak{S}(\{\techsiteTS: s\in\siteSet_r\})$ & $\times$ & Order of sites visited by tech $t$ of $r$ \\
    $\worktimeT$ & $\IntervalVar$ & $ [[0,\mwDur],[0,\mwDur]] $ & $\checkmark$ &  Working time of tech $t$ of $r$\\
    $\techShiftDurT$ & $\IntegerVar$ & $\duration$ & $\checkmark$ & Shift duration of tech $t$ of $r$\\
    $\numTS$ & $\IntegerVar$ & $\displaystyle[0,\sum_{s'}\overline{m}_{ss'w}]$ & $\checkmark$ & $\#$ of endpoints migrated by tech $t$ of $r$ at site $s$ \\
    $\numTSSp$ & $\IntegerVar$ & $[0,m_{ss'w}]$ & $\checkmark$ & $\#$ of endpoints migrated by tech $t$ of $r$ at site $s$ for  $\{s,s'\}$\\
\bottomrule
\end{tabular}    
\end{table}
Note that, although the variables are initially defined for all technicians, they are mostly ``optional'', meaning that if no shift is assigned to a technician $t$ from a region $r$, then its associated variables $\techsiteTS, \worktimeT, \techShiftDurT, \numTS $ and $\numTSSp$ are absent from the solution and vice-versa. For interval variables $\techsiteTS, \worktimeT$, the initial domain of their ``start'' and ``end'' times are $[0,\mwDur]$, the maximum time interval possible for a maintenance window. The initial domain of $\seqT$ is  $\mathfrak{S}(\{\techsiteTS: s\in\siteSet_r\})$, the set of permutations from the interval variables $\techsiteTS$. 
With the objective of minimizing the cost of a plan, the model $\LBBDSP = \CP$ is:\\[-8mm]
\begin{subequations}
\label{eq:CP_model}
\begin{align}
\min \ \  &\cost\sum_{r\in\regSet}\sum_{t\in\T_{rw}}\techShiftDurT  \label{eq:CP-obj}\\
\text{s.t.}\ \ 
&\NoOverlap(\seqT, \travel)&&t\in\T_{rw}, r\in\regSet\label{eq:no_overlap}\\
& \Span(\worktimeT, \{\techsiteTS: s\in\siteSet_r\}) &&t\in\T_{rw}, r\in\regSet\label{eq:span}\\
& \IfThen(\PresenceOf(\worktimeT), \EndOf(\worktimeT) \leq  \techShiftDurT) && t\in\T_{rw}, r\in\regSet\label{eq:CP_shift_duration}\\
& \PresenceOf(\techsiteTS) = \PresenceOf(\numTS) && s\in\siteSet_r,t\in\T_{rw}, r\in\regSet\label{eq:CP_tech_at_s_2}\\
&\IfThen(\PresenceOf(\techsiteTS), \LengthOf(\techsiteTS) = \durmigr\numTS)  && s\in\siteSet_r,t\in\T_{rw}, r\in\regSet\label{eq:CP_tech_at_s_3}\\
&\sum_{s':\{s, s'\}\in\sitePair}\numTSSp=\numTS && s\in\siteSet_r,t\in\T_{rw}, r\in\regSet\label{eq:CP_nss_of_t}\\
&\numTSSp=\numTSpS && \{s, s'\}\in\sitePair,t\in\T_{rw}, r\in\regSet\label{eq:CP_just_one_s_of_pair}\\
&\sum_{r\in\regSet}\sum_{t\in\T_{rw}}\numTSSp = \overline{m}_{ss'w} && \{s, s'\}\in\sitePair\label{eq:CP_tech_m}\\
& \sum_{r\in\regSet}\sum_{t\in\T_{rw}} \PresenceOf(\worktimeT) \leq \aEng \maxEng.\label{eq:CP_eng_max}
\end{align}
\end{subequations}       
By default, the interval variables in a $\SequenceVar$ are allowed to have overlaps. Using constraints \eqref{eq:no_overlap}, we make sure that the time intervals a technician $t$ spends at different sites are disjoint, and moreover these interval are defined apart enough such that there is time for the technician to travel from one site to another (as given by the matrix of travel times  $T$). Constraints \eqref{eq:span} make the working time of the technicians consistent with the planned activities. For the working technicians in the current plan, constraints  \eqref{eq:CP_shift_duration} determine their shift duration from $\duration$. Constraints \eqref{eq:CP_tech_at_s_2} and \eqref{eq:CP_tech_at_s_3} together ensure that the circuit endpoints are only distributed among the working technicians and the time a technician spends at a site is the amount of time it needs for migrating the endpoints. Constraints \eqref{eq:CP_nss_of_t} determine the site pairs to which the migrated endpoints belong, while constraints \eqref{eq:CP_just_one_s_of_pair} prevent the technician to move both ends of the circuits for the same site pair. Via constraints \eqref{eq:CP_tech_m}, the number of planned circuit migrations between the site pairs, as given by the master problem, is respected in the CP solution. Through constraints \eqref{eq:CP_eng_max}, the number of technicians to be used in each plan is bounded above by the maximum number of available engineers. 

\subsection{LBBD Cuts}\label{sec:LBBD_cuts}
Depending on the status of the $\CP$, we might need LBBD feasibility/optimality cuts to proceed. Next, we design these cuts and discuss their merits and drawbacks.

\subsubsection{LBBD Feasibility Cuts.}\label{sec:feas-lbbd-cut}
Infeasibility of $\CP$ implies that there is not enough resources (technicians and/or engineers) to migrate $\overline{\boldsymbol{m}}_w$ circuit endpoints in $w$. Using LBBD feasibility cuts, we remove this solution from the master problem $\LBBDMP$, along with any other solution that requires at least the same amount of resources as $\overline{\boldsymbol{m}}_w$. $\LBBDMP$ solutions are general integer, so the usual type of (strengthened) no-good cuts are not applicable to our problem. As such, for each generated $\overline{m}_{ss'w}$ we create a new binary variable to compare the  solution $m_{ss'w}$ with $\overline{m}_{ss'w}$. At iteration $\iter$ of the LBBD algorithm, we add the following LBBD feasibility cuts to $\feascuts$:
\begin{subequations}
\label{eq:LBBD_feas}
\begin{align}
&m_{ss'w} - (\cirSS+1) \greater^{\iter}_{ss'w} \leq \overline{m}_{ss'w}-1, \quad \{s,s'\}\in\sitePair: \overline{m}_{ss'w}>0 \label{eq:greater}\\
&\sum_{\substack{\{s,s'\}\in\sitePair\\\overline{m}_{ss'w}>0}} \greater^{\iter}_{ss'w} \leq |\sitePair| - 1\\
&\boldsymbol{\greater}^{\iter}_w \in\{0,1\}^{|\sitePair|}\label{eq:greater-domain}
\end{align}
\end{subequations}
As the number of iterations grows, the number of binary variables created via \eqref{eq:LBBD_feas} and the additional constraints can grow prohibitively large. However, as seen in our numerical experiments, the fact that the $\LBBDMP$ solutions first go through the Benders decomposition phase (with its own feasibility checks), refines them and reduces the number of required feasibility cuts. Additionally, \eqref{eq:LBBD_feas} cut-off not just one solution, but a potentially large set of solutions, further reducing the chance of needing new feasibility cuts. Proof of Proposition \ref{thm:LBBD_feas} explains the logic behind the design of our LBBD feasibility cuts.
\begin{proposition}\label{thm:LBBD_feas}
The set of inequalities \eqref{eq:LBBD_feas} are valid LBBD feasibility cuts for $\LBBDMP$.
\end{proposition}
\proof{
We introduce a vector of binary variables $\greater_{ss'w}^{\iter}, \{s,s'\}\in\sitePair$ that takes 1 if $m_{ss'w}\geq \overline{m}_{ss'w}$. We have the following logical relationships:
\begin{align*}
    m_{ss'w}\geq \overline{m}_{ss'w} \Longrightarrow \greater_{ss'w}^{\iter} = 1 &\quad {\equiv}\quad && \greater_{ss'w}^{\iter}=0 \Longrightarrow m_{ss'w}< \overline{m}_{ss'w} & (\text{contraposition})\\
    &\quad {\equiv}\quad && \greater_{ss'w}^{\iter}=0 \Longrightarrow m_{ss'w}\leq \overline{m}_{ss'w}-1 & (m_{ss'w}\in\mathbb{Z}_+)\\
    &\quad {\equiv}\quad && m_{ss'w} -(\cirSS + 1) \greater^\iter_{ss'w} \leq \overline{m}_{ss'w}-1 & (m_{ss'w}\leq\cirSS)
\end{align*}
The last inequality becomes redundant if $\greater^\iter_{ss'w} = 1$. In order to cut-off any  solution $\boldsymbol{m}_w\geq\overline{\boldsymbol{m}}_w$, it is sufficient to make sure that the number of migrated circuit endpoints between at least one site pair $\{s,s'\}$ is less than  $\overline{m}_{ss'w}$, which completes the proof. 
}

\subsubsection{LBBD Optimality Cuts.} Assuming that $\CP$ is feasible and solved to optimality, let $\CPSol$ be its optimal objective value. If $\overline{\eta}_w \geq \CPSol$ for a $w\in\wSet$, then $(\overline{\boldsymbol{m}}, \overline{\eta}_w)$ is feasible for $\LBBDMP$. If $\overline{\eta}_w \geq \CPSol$ for all $w\in\wSet$, then $\overline{\boldsymbol{m}}$ is optimal for the problem and the algorithm stops. Otherwise, for any $w$ with $\overline{\eta}_w < \CPSol$, we cut-off the solution pair ($\overline{\eta}_w , \overline{\boldsymbol{m}}_w$). At iteration $\iter$, we add the following LBBD optimality cuts to $\optcuts$:
\begin{subequations}
\label{eq:LBBD_opt}
\begin{align}
&\eqref{eq:greater}, \eqref{eq:greater-domain}\\
&\displaystyle\eta_w \geq \CPSol \Big(1- \sum_{\substack{\{s,s'\}\in\sitePair\\\overline{m}_{ss'w}>0}} \big(1-\greater^{\iter}_{ss'w}\big)\Big)
\end{align}
\end{subequations}
As with the feasibility cuts \eqref{eq:LBBD_feas}, we are creating a large number of new binary variables and constraints. Once more, our numerical analysis show that, for the considered instances, the number of iterations required for the algorithm to converge to a reasonable gap is small enough  that $\LBBDMP$ remains scalable despite the large number of binary decision variables and constraints added at each iteration.
\begin{proposition}
The set of inequalities \eqref{eq:LBBD_opt} are valid LBBD optimality cuts for $\LBBDMP$.
\end{proposition}
\proof{
The LBBD optimality cut is valid if ($i$) it can determine whether the solution $\overline{\boldsymbol{m}}_w$ is optimal, ($ii$) returns a valid bound for the solutions other than $\overline{\boldsymbol{m}}_w$, and ($iii$) it is tight at $\overline{\boldsymbol{m}}_w$. Because the $\text{CP}_w(\boldsymbol{m})$'s cost function is non-decreasing, for any $\boldsymbol{m}_w\geq\overline{\boldsymbol{m}}_w$  we have $\overline{\text{CP}}_w(\boldsymbol{m})\geq\CPSol$. With the binary variables $\greater^\iter_{ss'w}$ as defined in the proof of Proposition \ref{thm:LBBD_feas},  it is clear that inequalities \eqref{eq:LBBD_opt} are tight at $\overline{\boldsymbol{m}}_w$. It also imposes a negative bound on $\boldsymbol{m}_w < \overline{\boldsymbol{m}}_w$, and $\CPSol$ if $\boldsymbol{m}_w > \overline{\boldsymbol{m}}_w$, which completes the proof.
}

We can also characterize a set of solutions other than $\overline{\boldsymbol{m}}_{w}$ for which the cuts \eqref{eq:LBBD_opt} are tight: If $\hat{m}_{ss'w}-\overline{m}_{ss'w}>0$, then at least $\theta (\hat{m}_{ss'w}-\overline{m}_{ss'w})$ minutes should be added to $\worktimeT$ for some technician $t$. Let $\overline{n}^\textsc{s}_{rts}, \overline{\textsc{wtime}}_{rt}, \overline{\Delta}_{\textsc{shift}}^{rt}$ be the solutions from the CP model. Define a bipartite graph $\mathcal{G} = (\mathcal{V}_1\cup\mathcal{V}_2,\mathcal{E})$, with the following set of nodes:
$$\mathcal{V}_1 = \bigcup_{r\in\regSet}\T_{rw}\setminus\{t\in\T:\ \overline{n}^\textsc{s}_{rts}=0\},\quad \mathcal{V}_2 = \{\{s,s'\}\in\sitePair:\ \hat{m}_{ss'w}-\overline{m}_{ss'w}>0\}.$$
In other words, $\mathcal{V}_1$ is the set of working technicians and $\mathcal{V}_2$ is the set of site pairs that need extra circuit migration in the new solution. We build the set of links as follows: 
$$\mathcal{E} = \big\{(t,\{s,s'\}) \in \mathcal{V}_1\times\mathcal{V}_2:\ \exists  t\in\T_{\siteSet_rw},\ ,\  \overline{\textsc{wtime}}_{rt}+\theta (\hat{m}_{ss'w}-\overline{m}_{ss'w}) \leq \overline{\Delta}_{\textsc{shift}}^{rt} \big\},$$
i.e, we create  a link between a technician that is already planned to work at $s$ and has enough spare time within its shift to migrate the remaining circuits. If there exists a maximum matching of size $|\mathcal{V}_2|$, then $\CPSol$ does not change and \eqref{eq:LBBD_opt} is tight at $\hat{\boldsymbol{m}}_{w}$. Note that, the proposition is not biconditional.

\begin{figure}[t]
	\centering
	\scalebox{0.75}{
	\input{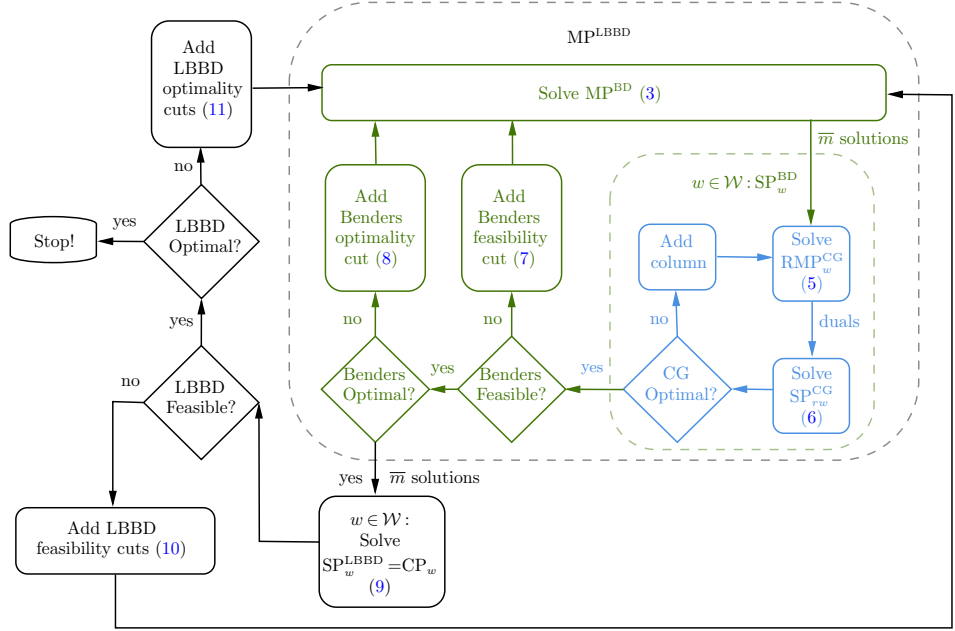}
	}
	\caption{Overall framework of the LBBD method for the NMP.}
	\label{fig:alg}
\end{figure}

Figure \ref{fig:alg}  presents the overall solution framework. Our LBBD method converges with an exact solution, or returns ``infeasible'' if the original problem is infeasible. This follows from the integrality of $\boldsymbol{m}$, which makes the set of feasible solutions for $\LBBDMP$ finite. Using the feasibility and optimality cuts \eqref{eq:LBBD_feas} and \eqref{eq:LBBD_opt}, we either cut-off the solution passed by $\LBBDMP$ or stop the algorithm if no action is needed after checking the feasibility and optimality conditions. When the algorithm stops, if we do not have any incumbent $\overline{\boldsymbol{m}}$, the original problem is infeasible, otherwise the algorithm has reached to optimality. 

\section{Numerical Results}
\label{sec:numer_results}

In this section we have evaluated our LBBD method with the NMP instances that are built over six networks of different sizes, with some deployed over a continent and another serving a region.
In Section \ref{sec:datasets}, we provide some details on the characteristics of the instances. This is followed by algorithmic and cost analysis discussions in Sections \ref{sec:alg_analysis} and  \ref{sec:cost_analysis}, respectively.

\subsection{Data Sets}\label{sec:datasets}
\begin{figure}[hbt]
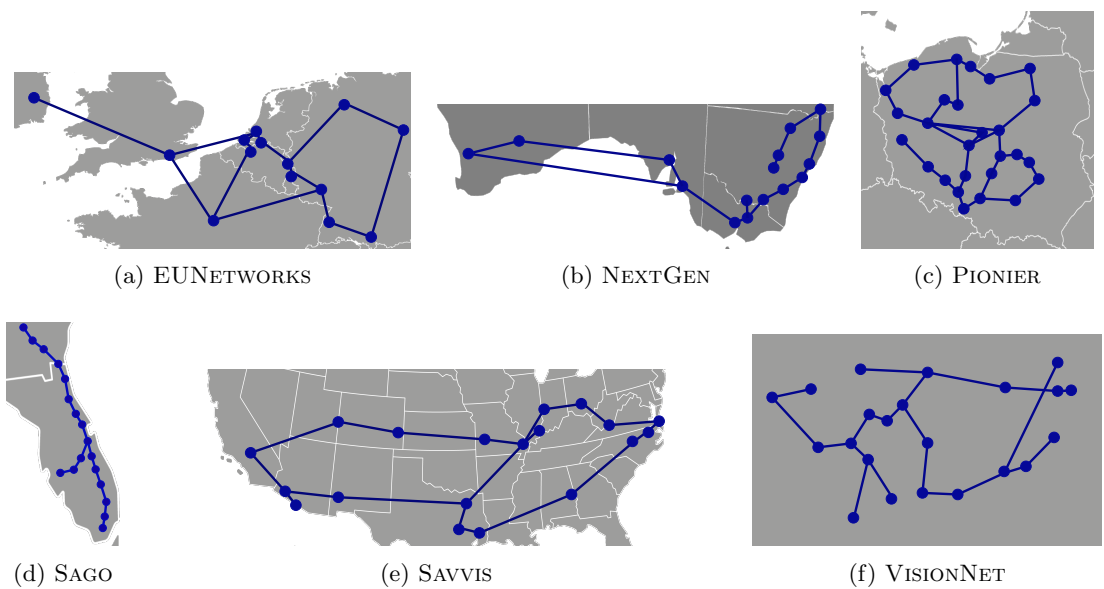

	\centering
		\subfloat[\EUnet{}]{\scalebox{0.4}{\input{Figures/EU}}} \ 
	    \subfloat[\NextGen{}]{\scalebox{0.4}{\input{Figures/NextGen}}}\
        \subfloat[\Pionier{}]{\scalebox{0.4}{\input{Figures/Pionier}}}\\
		\subfloat[\Sago{}]{\scalebox{0.4}{\input{Figures/Sago}}} \hspace{8mm} 
		\subfloat[\Savvis{}]{\scalebox{0.4}{\input{Figures/Savvis}}}\hspace{8mm} 
        \subfloat[\Vision{}]{\scalebox{0.4}{\input{Figures/VisionNet}}}\\[2mm]
		\caption{Topology of the base networks}
		\label{fig:networks}
\end{figure}

In Figure \ref{fig:networks} we have reproduced the topology of six real optical fiber networks as reported by \cite{6027859}. Considering the nodes of the graphs as the sites of the network, we  use the nodes' coordinates to partition the sites into regions, using the quality threshold clustering algorithm with 80km threshold. Detailed characteristics of the base networks are given in Table \ref{tab:networks}, including the number of sites, regions, maximum number of sites in a region, their geographic extent and location.  
\begin{table}[htb]
  \centering
  \caption{Characteristics of base networks}
  \scalebox{0.9}{\renewcommand{\arraystretch}{1.25}
    \small
    \begin{tabular}{p{3cm}cccccc}
    \toprule
    Characteristics & \EUnet & \NextGen & \Pionier & \Sago & \Savvis & \Vision  \\
    \midrule
    $|\siteSet|\ $ & 15  & 16    & 21  & 18    & 19    & 22     \\
    $|\regSet|\ $ & 11  & 15    & 16  & 10    & 19    & 16     \\
    $|\siteSet_r|^{\max}\ $ &  4  & 2     & 3    &    3   & 1     &   3    \\
    Geographic scope & Continent & Country & Country & Region & Country & Region\\
    Coverage & Western Europe & Australia & Poland & Florida \& Georgia & USA & Montana\\
    \bottomrule
    \end{tabular}%
    }
  \label{tab:networks}%
\end{table}%
Let $\mu_{EP/s}$ and $\sigma_{EP/s}$ be the mean and standard deviation of the number of circuit endpoints at a site, respectively. Using a Lognormal distribution with parameters $\mu_{EP/s}$ and $\sigma_{EP/s}$, we have generated 114 NMP instances over the base networks by creating six instances per network and  $(\mu_{EP/s},\sigma_{EP/s})$. The considered values for $\mu_{EP/s}$ are comparable with $\mu_{EP/s} = 6.18$ of a real data set studied in \cite{pouya2017minimum} (Ciena data set I). The instances then are characterized by their base network, the number of circuit endpoints and site pairs with at least one circuit.  Table \ref{tab:instances} presents a summary of the characteristics of the instances, and the details are provided in the e-companion.
\begin{table}[htb]
  \centering
  \caption{Characteristics of the NMP instances}
    \scalebox{0.9}{\renewcommand{\arraystretch}{1.25}
    \small
\begin{tabular}{lcccccrrcc}
\toprule
\multirow{2}[2]{*}{Dataset} & \multirow{2}[2]{*}{\#Instances} & \multirow{2}[2]{*}{$|\wSet|$} & \multirow{2}[2]{*}{$ \maxCir $} & \multirow{2}[2]{*}{$\mu_{EP/s}$ } & \multirow{2}[2]{*}{$\sigma_{EP/s}$ } & \multicolumn{2}{c}{Endpoints} & \multicolumn{2}{c}{$| \sitePair| $} \\
\cmidrule(lr){7-8}\cmidrule(lr){9-10}      &       &       &       &       &       & $\min$   & $\max$   & $\min$   & $\max$ \\
\midrule
\EUnet5 & 6     & 3     & 30    & 5     & 2.5   & 46    & 66    & 22    & 32 \\
\EUnet6 & 6     & 3     & 30    & 6     & 3.0   & 70    & 88    & 28    & 34 \\
\EUnet7 & 6     & 4     & 40    & 7     & 3.5   & 80    & 98    & 34    & 40 \\
\EUnet8 & 6     & 4     & 40    & 8     & 4.0   & 94    & 116   & 38    & 43 \\
\EUnet9 & 6     & 4     & 40    & 9     & 4.5   & 104   & 150   & 41    & 46 \\
\midrule
\NextGen5 & 6     & 3     & 30    & 5     & 2.5   & 52    & 74    & 23    & 31 \\
\NextGen6 & 6     & 3     & 30    & 6     & 3.0   & 62    & 96    & 26    & 39 \\
\NextGen7 & 6     & 3     & 30    & 7     & 3.5   & 78    & 92    & 31    & 37 \\
\NextGen8 & 6     & 4     & 40    & 8     & 4.0   & 92    & 126   & 35    & 47 \\
\NextGen9 & 6     & 5     & 40    & 9     & 4.5   & 98    & 134   & 41    & 51 \\
\midrule
\Pionier5 & 6     & 3     & 30    & 5     & 2.5   & 70    & 98    & 35    & 43 \\
\Pionier6 & 6     & 3     & 30    & 6     & 3.0   & 96    & 112   & 40    & 51 \\
\midrule
\Sago5 & 6     & 3     & 30    & 5     & 2.5   & 64    & 90    & 29    & 39 \\
\Sago6 & 6     & 3     & 40    & 6     & 3.0   & 80    & 94    & 32    & 40 \\
\midrule
\Savvis5 & 6     & 3     & 30    & 5     & 2.5   & 54    & 102   & 25    & 42 \\
\Savvis6 & 6     & 3     & 30    & 6     & 3.0   & 84    & 106   & 38    & 41 \\
\Savvis7 & 6     & 3     & 30    & 7     & 3.5   & 88    & 120   & 37    & 47 \\
\midrule
\VisionNet5 & 6     & 3     & 30    & 5     & 2.5   & 82    & 104   & 36    & 47 \\
\VisionNet6 & 6     & 4     & 40    & 6     & 3.0   & 106   & 126   & 47    & 52 \\
\bottomrule
\end{tabular}%
    }
  \label{tab:instances}%
\end{table}%
For the travel time between the sites, we consider an average speed of 80km/h. The migration duration of a single endpoint is 20 minutes. Technicians and engineers are paid 108\$ and 140\$ per hour, respectively. Technicians may work in shifts of $\durshift\in\{6,8\}$ hours and an engineer can coordinate 5 technicians at a time. These instance parameters are based on the inputs of a real data set analyzed in \cite{pouya2018new}.

\subsection{Implementation Details}\label{sec:imp_detail}
Programs are written in Python and run on Niagara supercomputer servers   \citep{loken2010scinet, ponce2019deploying} using Gurobi 9.1.1 as the MIP solver \citep{gurobi} and CP Optimizer 20.1 as the CP solver \citep{cplex}.  For implementation of the Benders decomposition, we have used Gurobi's callback feature such that, while solving the $\BDMP$, at every integer node of the branch-and-bound tree, we separate the Benders optimality and feasibility cuts by solving the CG subproblems. In implementing the LBBD method, we have opted for a cutting plane framework where every time after the convergence of the Benders decomposition method, we solve a CP problem to separate the violated LBBD optimality and feasibility cuts, if any. The reason for employing the cutting plane framework in this case is that our LBBD feasibility and optimality cuts introduce new binary variables. This is not allowed in the Gurobi's callback implementation, since introducing new integer variables changes the structure of the branch-and-bound tree.

The programs stop whenever the algorithm reaches an optimality gap under 10\%, or the running time passes 3 hours, whichever happens first. Our preliminary experiments showed that at the beginning of the solution process, early stopping of the Benders decomposition step reduces the computational effort, as it provides the $\LBBDMP$ with informative cuts even without having to solve it to optimality. Therefore, in solving the $\LBBDMP$ (i.e., the Benders decomposition phase) we start with a MIP gap of 10\%, and reduce this by 5\% at each iteration. Note that in these steps we should use the best lower bound on the solution in measuring the gap, which gives us a valid lower bound on the NMP instances. 
We have also made the following implementation choices and analyzed their impacts in the next section:
\begin{itemize}
    \item In solving the CG subproblems, we use the hybrid CG as explained in Remark \ref{rem:hybriddCG}.
    \item Despite the fact that the CG subproblems are frequently solved  during the algorithm, we do not remove the generated columns from $\Gp$.
    \item Considering the symmetry between the maintenance windows, we propagate the Benders optimality cuts, i.e., every time a Benders optimality cut is separated for a certain maintenance window $w$, we also modify and apply it for all $w'\in\wSet\setminus\{w\}$. 
\end{itemize}

\subsection{Algorithmic Analysis}
\label{sec:alg_analysis}
In the following sections, we first examine our implementation choices, and based on the best obtained settings, solve the NMP instances. We break down the algorithmic analysis of the method into each solution paradigm and study their share of the computational effort.  
\subsubsection{Performance of the Hybrid CG.} \label{sec:hybrid-analysis}
We have solved the first instances of \EUnet5, \NextGen5, \Pionier5, \Sago5, \Savvis5 and \Vision5 with a pure CP-based CG, as well as the hybrid version where first a MIP is solved to generated improving columns with ordered paths, and then the CP pricing problem \eqref{eq:pp} continues the process by generating columns with general paths or verifies the optimality. The results in Table \ref{tab:two-sub-perform} show that  the hybrid CG reduces the computational effort of solving the $\LBBDMP$ by orders of magnitude. Most significantly, in 3 hours the pure CP-based method is unable to achieve  a gap below 10\% for \Pionier5, \Sago5 and \Vision5, while the hybrid CG  in most cases closes  the gap in a  short amount of time. As mentioned in Section \ref{sec:CG}, in our backbone networks many regions have one or two sites, so the MIP pricing problem is already sufficient to guarantee their optimality. This greatly benefits the geographically larger networks \EUnet, \NextGen{} and \Savvis{} to the point that the instances based on the last two networks can be solved to optimality by only solving the MIP subproblem.
\begin{table}[htbp]
  \centering
  \caption{Impact of the hybrid CG on the overall performance of the method}
  \scalebox{0.9}{\renewcommand{\arraystretch}{1.25}
  \small
\begin{tabular}{lrrrcrrrc}
\toprule
\multirow{3}[3]{*}{Dataset} & \multicolumn{4}{c}{CP-sub}    & \multicolumn{4}{c}{CP/MIP-sub} \\
\cmidrule(lr){2-5}\cmidrule(lr){6-9}      & \multirow{2}[2]{*}{\# Columns} & \multicolumn{2}{c}{Time (s)} & \multirow{2}[2]{*}{Gap (\%)} & \multirow{2}[2]{*}{\# Columns} & \multicolumn{2}{c}{Time (s)} & \multirow{2}[2]{*}{Gap (\%)} \\
\cmidrule(rl){3-4}\cmidrule(rl){7-8}      &       & \multicolumn{1}{c}{CG} & \multicolumn{1}{c}{$\LBBDMP$} &       &       & \multicolumn{1}{c}{CG} & \multicolumn{1}{c}{$\LBBDMP$} &  \\
\midrule
\EUnet5 & 138   & 4184.5 & 4192.2 & 0.0\% & 509   & 627.1 & 643.0 & 0.0\% \\
\NextGen5 & 45    & 304.2 & 305.8 & 0.0\% & 45    & 5.2   & 10.4  & 0.0\% \\
\Pionier5 & 493   & \texttt{Timeout} & \texttt{Timeout} & 47.4\% & 993   & 1039.4 & 1181.1 & 0.8\% \\
\Sago5 & 314   & \texttt{Timeout} & \texttt{Timeout} & 31.6\% & 141   & 655.2 & 692.5 & 0.0\% \\
\Savvis5 & 57    & 1219.7 & 1257.3 & 0.0\% & 57    & 11.7  & 78.2  & 0.0\% \\
\VisionNet5 & 512   & \texttt{Timeout} & \texttt{Timeout} & 50.0\% & 1738  & 1878.2 & 2010.4 & 3.0\% \\
\bottomrule
\end{tabular}%
    }
  \label{tab:two-sub-perform}%
\end{table}

The comparatively remarkable performance of the hybrid CG can further be explained by a common practice in the column generation literature, where sub-optimal columns are also added to the restricted master problem in order to make it more descriptive and obtain more informative dual solutions for the pricing problem. From Table \ref{tab:two-sub-perform} we clearly see that the number of generated columns in the hybrid CG is generally much larger, meaning that the MIP subproblem is returning many columns with sub-optimal reduced-costs, and we are practically following the mentioned strategy in increasing the number of columns. 

\subsubsection{CG Column Management}
The CG subproblems are solved frequently, at every integer node of the branch-and-bound tree every time the algorithm reaches the Benders decomposition stage. Furthermore, the models are modified constantly, with different right-hands sides in each node.  It is not immediately clear whether keeping all the columns leads to a speedup of the method or overloading the $\CGMP$. For the same instances of Section \ref{sec:hybrid-analysis}, we have evaluated the impact of removing the columns after processing an integer node of the tree.
Results in Table \ref{tab:cg-col-manage} show that the overall solution process indeed benefits from previously generated columns, and we are still able to efficiently solve the $\CGMP$. 
\begin{table}[b]
  \centering
  \caption{Performance of the LBBD method under two column management strategies}
  \scalebox{0.9}{\renewcommand{\arraystretch}{1.25}
  \small
\begin{tabular}{lrrrrrr}
\toprule
\multirow{3}[3]{*}{Dataset} & \multicolumn{3}{c}{CG w/o previous columns} & \multicolumn{3}{c}{CG w/ all generated columns} \\
\cmidrule(lr){2-4}\cmidrule(lr){5-7}      & \multirow{2}[2]{*}{\ \#\textsc{iter.}} & \multicolumn{2}{c}{Time (s)} & \multirow{2}[2]{*}{\ \#\textsc{iter.}} & \multicolumn{2}{c}{Time (s)} \\
\cmidrule(rl){3-4}\cmidrule(rl){6-7}      &       & \multicolumn{1}{c}{CG} & \multicolumn{1}{c}{$\LBBDMP$} &       & \multicolumn{1}{c}{CG} & \multicolumn{1}{c}{$\LBBDMP$} \\
\midrule
\EUnet5 & 4     & 727.3 & 1294.1 & 4     & 627.1 & 643.0 \\
\NextGen5 & 2     & 307.6 & 308.8 & 14    & 5.2   & 10.4 \\
\Pionier5 & 6     & 1162.1 & 2132.9 & 5     & 1039.4 & 1181.1 \\
\Sago5 & 10    & 1601.3 & 2134.5 & 10    & 655.2 & 692.5 \\
\Savvis5 & 3     & 74.0  & 139.8 & 4     & 11.7  & 78.2 \\
\VisionNet5 & 14    & 2572.0 & 7933.0 & 7     & 1878.2 & 2010.4 \\
\bottomrule
\end{tabular}%
    }
  \label{tab:cg-col-manage}%
\end{table}%
However, the results also show that solving the CG models with an initial set $\Gp = \emptyset$ is still manageable and if the size of the instances leads to poor performance of the $\CGMP$, we can switch to removing the columns or design a hybrid column management strategy.

\subsubsection{Propagation of the Optimality Cuts.} \label{sec:propagation}
The symmetry among the maintenance windows $w\in\wSet$ in the NMP signifies that if an inequality $\eta_w\geq \boldsymbol{g}(\overline{\eta}_w , \overline{\boldsymbol{m}}_w)\boldsymbol{m}_w$ cuts  off the solution pair ($\overline{\eta}_w , \overline{\boldsymbol{m}}_w$) from the solution space, a similar inequality $\eta_{w'}\geq \boldsymbol{g}(\hat{\eta}_{w'} , \hat{\boldsymbol{m}}_{w'})\boldsymbol{m}_{w'}$ removes $(\hat{\eta}_{w'} , \hat{\boldsymbol{m}}_{w'})$  if $\hat{\eta}_{w'}=\overline{\eta}_{w}, \hat{\boldsymbol{m}}_{w'}=\overline{\boldsymbol{m}}_{w}$. Therefore, any optimality cut obtained for one $w$ can be modified to use for all maintenance windows. This is called Benders cut propagation \citep{roshanaei2017propagating}. In Table \ref{tab:propagate_cuts} we have examined the effect of propagating the Benders optimality cuts in solving the instances of Section \ref{sec:hybrid-analysis} and compared the results with a version of the method that only uses the classic local Benders optimality cuts. Columns \#Cut$^{\textsc{BD}}_{\textsc{opt}}$ and \#Cut$^{\textsc{LBBD}}_{\textsc{opt}}$ respectively report the number of the Benders and LBBD optimality cuts at the end of the solution process.
\begin{table}[t]
  \centering
  \caption{Performance of the LBBD method with and without propagation of the Benders optimality cuts}
  \scalebox{0.9}{\renewcommand{\arraystretch}{1.25}
  \small
\begin{tabular}{lrrrrrrrr}
\toprule
\multirow{2}[2]{*}{Dataset} & \multicolumn{4}{c}{LBBD with local Benders cuts} & \multicolumn{4}{c}{LBBD with propagated Benders cuts} \\
\cmidrule(lr){2-5}\cmidrule(lr){6-9}      & \multicolumn{1}{c}{\ \#\textsc{iter.}\ } & \multicolumn{1}{c}{\ \#Cut$^{\textsc{BD}}_{\textsc{opt}}$} & \multicolumn{1}{c}{\ \#Cut$^{\textsc{LBBD}}_{\textsc{opt}}$} & \multicolumn{1}{c}{Time (s)} & \multicolumn{1}{c}{\#\textsc{iter.}\ } & \multicolumn{1}{c}{\ \#Cut$^{\textsc{BD}}_{\textsc{opt}}$} & \multicolumn{1}{c}{\ \#Cut$^{\textsc{LBBD}}_{\textsc{opt}}$} & \multicolumn{1}{c}{Time (s)} \\
\midrule
\EUnet5 & 3     & 49    & 81    & 624.9 & 4     & 42    & 139   & 643.0 \\
\NextGen5 & 4     & 16    & 170   & 8.8   & 14    & 32    & 579   & 10.4 \\
\Pionier5 & 9     & 241   & 567   & 1841.3 & 5     & 90    & 269   & 1181.1 \\
\Sago5 & 22    & 726   & 1297  & 3507.9 & 10    & 91    & 553   & 692.5 \\
\Savvis5 & 3     & 147   & 101   & 71.8  & 4     & 45    & 174   & 78.2 \\
\VisionNet5 & 13    & 723   & 836   & 5835.2 & 7     & 202   & 359   & 2010.4 \\
\bottomrule
\end{tabular}%
    }
  \label{tab:propagate_cuts}%
\end{table}%
The results indicate that, the cut propagation is quite effective in accelerating the solution process of the instances that are based on geographically smaller networks \Pionier, \Sago{} and \VisionNet, while its performance is comparable to the classic version for the networks stretched over a larger geographic scope  \EUnet, \NextGen{} and \Savvis. Note that, the same arguments apply to the LBBD optimality cuts, yet any gain from propagating the cuts \eqref{eq:LBBD_opt} is outweighed by their numbers. As a result,  we only propagate the Benders optimality cuts in the LBBD method.

\subsubsection{Computational Performance.} 
\begin{figure}[b]
	{\footnotesize\scalebox{0.95}{
	\input{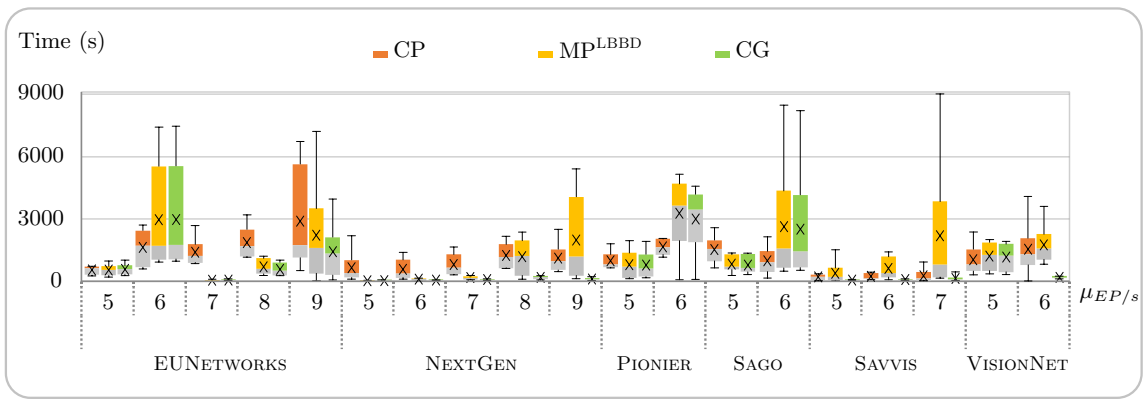}}}
	\caption{Breakdown of the solution times spent at each step of the decomposition. The time to solve the $\LBBDMP$ includes the total solution time of the CG models.}
	\label{fig:time_analysis}
\end{figure}

With the implementation choices fixed by the analyses of Sections  \ref{sec:hybrid-analysis} - \ref{sec:propagation}, we have solved the NMP instances of Table \ref{tab:instances} and reported the confidence intervals on the optimality gaps, solution costs, the number of LBBD iterations and solution times in Table \ref{tab:sol-time}. 
\begin{table}[htb]
  \centering
  \caption{Solution and computational effort}
  \scalebox{0.9}
      {\renewcommand{\arraystretch}{1.2}
  \small
    \begin{tabular}{lcccrrcrr}
    \toprule
    \multirow{2}[4]{*}{Dataset} & \multicolumn{2}{c}{Gap} & \multicolumn{2}{c}{Cost (\$)} & \multicolumn{2}{c}{\#\textsc{iter.}} & \multicolumn{2}{c}{Time (s)} \\
\cmidrule(lr){2-3} \cmidrule(lr){4-5}\cmidrule(lr){6-7} \cmidrule(lr){8-9}            & Mean  & CI width & Mean  & CI width & Mean  & CI width & Mean  & CI width \\
    \midrule
    \EUnet5 & 0.6\% & 1.1\% & 10472.0 & 448.7 & 3.5   & 0.6   & 1124.1 & 268.6 \\
    \EUnet6 & 1.0\% & 1.9\% & 11424.0 & 923.4 & 10.5  & 4.3   & 4607.1 & 2309.7 \\
    \EUnet7 & 4.7\% & 1.9\% & 11877.3 & 648.9 & 6.5   & 2.0   & 1459.0 & 509.7 \\
    \EUnet8 & 4.9\% & 2.9\% & 12693.3 & 611.3 & 7.3   & 1.6   & 2582.6 & 539.9 \\
    \EUnet9 & 3.7\% & 3.0\% & 13328.0 & 615.6 & 10.8  & 6.1   & 5109.8 & 3376.7 \\
    \midrule
    \NextGen5 & 1.2\% & 2.2\% & 13192.0 & 954.9 & 6.0   & 3.0   & 699.2 & 560.3 \\
    \NextGen6 & 2.6\% & 2.1\% & 14960.0 & 1443.7 & 4.8   & 2.3   & 675.7 & 374.2 \\
    \NextGen7 & 2.4\% & 1.7\% & 15640.0 & 243.3 & 5.3   & 2.0   & 997.9 & 415.8 \\
    \NextGen8 & 1.8\% & 1.3\% & 16048.0 & 615.6 & 8.2   & 4.1   & 2445.6 & 399.6 \\
    \NextGen9 & 4.1\% & 2.5\% & 16864.0 & 721.8 & 8.0   & 2.9   & 3120.4 & 1939.1 \\
    \midrule
    \Pionier5 & 1.9\% & 1.7\% & 15912.0 & 821.6 & 5.7   & 1.6   & 1817.2 & 796.4 \\
    \Pionier6 & 3.2\% & 2.1\% & 17136.0 & 533.1 & 9.2   & 2.0   & 4929.4 & 1385.8 \\
    \midrule
    \Sago5 & 1.3\% & 1.5\% & 11832.0 & 625.1 & 8.3   & 2.3   & 2387.6 & 697.1 \\
    \Sago6 & 5.2\% & 2.4\% & 12648.0 & 498.7 & 5.7   & 2.3   & 3635.1 & 2640.4 \\
    \midrule
    \Savvis5 & 1.3\% & 1.4\% & 18224.0 & 1109.8 & 4.5   & 1.4   & 619.6 & 475.2 \\
    \Savvis6 & 3.1\% & 2.1\% & 19312.0 & 721.8 & 6.0   & 1.6   & 855.6 & 460.3 \\
    \Savvis7 & 4.2\% & 2.2\% & 19856.0 & 721.8 & 6.0   & 3.7   & 2472.6 & 2731.3 \\
    \midrule
    \VisionNet5 & 4.9\% & 1.8\% & 17272.0 & 877.3 & 5.8   & 3.5   & 2259.2 & 899.9 \\
    \VisionNet6 & 3.8\% & 1.9\% & 17952.0 & 377.0 & 7.3   & 3.1   & 3308.0 & 769.4 \\
    \bottomrule
    \end{tabular}%
    }
  \label{tab:sol-time}%
\end{table}%
In all instances, the LBBD method returns quality solutions with small gaps in a reasonable amount of time. The small number of iterations in the LBBD stage is a result of multiple levels of decomposition in the LBBD master problem, giving it enough information to pass on quality first-stage solutions to the subproblems. This is quite valuable, considering the expensive design of the LBBD feasibility and optimality cuts. In Figure \ref{fig:time_analysis}, we have further inspected the solution times. The box plots in this figure represent the distribution of the time that each step of the algorithm spends for solving the six instances of each dataset. We observe that in the smaller instances the $\LBBDMP$ run time is essentially dominated by the effort to solve the CG models. However, in larger instances with more endpoints, the time spent in solving the $\BDMP$ apart from the CG models becomes increasingly significant. In contrast, as the number of the endpoints increases, with the exception of \EUnet9, the computational effort by the CP models remains relatively stable with a mild increase. The exponential growth in the solution time of the $\LBBDMP$ is expected, considering the need to add more feasibility and optimality cuts, as reported in the e-companion.
It is worth noting that, none of the steps of the LBBD method dominates the others in all the studied instances. Accordingly, all solution paradigms play a determining role in the solution process and  enhancements in any of them can lead to overall improvement of the method.

\subsection{Network Migration Cost Analysis}
\label{sec:cost_analysis}
Figure \ref{fig:cost_analysis} demonstrates the impact of various parameters on the solutions of the \EUnet5  (similar analyses for two other networks are provided in the e-companion). 
\begin{figure}[b]
	\centering 
	\subfloat[\# Shifts vs. cost\label{fig:shift_cost}]{\scalebox{0.57}{
	\input{Figures/EU_ShiftCost}}}
	% -------------------
	\subfloat[Distribution of shifts\label{fig:shift_dur}]{\scalebox{0.57}{
	\input{Figures/EU_Shifts}
	}} \\
    % -------------------
	\subfloat[Cost vs. efficient working time\label{fig:working_cost}]{\scalebox{0.57}{
	\input{Figures/EU_WorkingCost}
	}}\\
	% -------------------
	\caption{Cost analysis for \EUnet}
	\label{fig:cost_analysis}
\end{figure}
As shown in Figure \ref{fig:shift_cost}, although the objective function has not benefited from engaging more technicians, increasing $\maxCir$ from 30 to 40 lowers the number of the shifts and migration costs. Finding the right amount of workforce is one of the main concerns of the network operators as it comes with logistics and hidden considerations. While a solution with higher cost does not necessarily imply having more shifts, fewer shifts might be considered an advantage as it requires fewer personnel and less human traffic in the sites. The distribution of the shifts, illustrated in Figure \ref{fig:shift_dur}, shows that the majority of the technicians are assigned shorter than 2-hour shifts if possible. This is the result of either shortage of technicians or circuits in the targeted regions. Since offering more technicians has not improved the migration costs, we can conclude that the large number of the 2-hour shifts is due to the availability of the circuit endpoints in the region of the working technicians. Once the 2-hour shifts are not offered, 4-hour shifts become the obvious choice. By moving to 4-hour shifts, the number of the shifts, as expected, declines and we deal with fewer number of the technicians during the migration. However, as mentioned before, this does not mean a decrease in the migration costs since the shifts are the indicators of the payments, not the working time. Figure \ref{fig:working_cost} displays the average working time (migrating circuit endpoints or traveling between the sites) of the technicians during the shifts. We notice that, by setting the minimum shift duration to 4 hours, the working time percentage is between 40\% - 50\% of a shift, while this number is above 60\% for the settings with minimum shift duration of 2 hours. The increase of the costs show that 4-hour shifts are on average less busy than 2-hour shifts. It should be mentioned that there are a few 6-hour shifts as well. Considering the negative impact of longer maintenance windows on the Service-Level-Agreement (SLA) between the network operators and the customers, the analysis indicates that having 6-hour shifts does not help the migration process in neither costs nor the number of the shifts. As seen in Figure \ref{fig:cost_analysis}, we can build equivalent solutions with shorter shifts and ignore the 6-hour ones.

\section{Concluding Remarks}\label{sec:conclusion}
Migrating legacy telecommunication networks to the latest technology involves planning synchronized technicians. We propose the first exact method for the network migration problem, a logic-based Benders decomposition method augmented by  column generation and constraint programming models. We also make several algorithmic enhancements and considerably improve the basic version of the algorithm with the classical adaption of Benders decomposition and column generation. The proposed solution method can further be adapted to any vehicle routing problem with node synchronization constraints. Our evaluations on the instances generated over six real networks show that our method is effective in obtaining quality solutions, and that all the solution paradigms are contributing to the efficiency of the method. Future work includes improving the scalability of the algorithm for the networks with a larger number of circuit endpoints, particularly the design of less expensive LBBD optimality and feasibility cuts. Furthermore, new methodological developments are needed to build migration plans in the presence of uncertainty, which is an inherent feature of such planning problems.

\section{Acknowledgment}
\small 
Computations were performed on the Niagara supercomputer at the SciNet HPC Consortium. SciNet is funded by: the Canada Foundation for Innovation; the Government of Ontario; Ontario Research Fund - Research Excellence; and the University of Toronto.

\bibliographystyle{informs2014}
\bibliography{Bibliography/Network_Migration,Bibliography/SVRP,Bibliography/b_p.bib,Bibliography/CP.bib}

\begin{thebibliography}{46}
\providecommand{\natexlab}[1]{#1}
\providecommand{\url}[1]{\texttt{#1}}
\providecommand{\urlprefix}{URL }

\bibitem[{Almughaless \protect\BIBand{} Alsaih(2010)}]{almughaless2010optimum}
Almughaless AA, Alsaih AM (2010) Optimum migration scenario from {PSTN} to
  {NGN}. \emph{IEEE International Conference on Communication Systems, Networks
  and Applications - ICCSNA}, volume~1, 227--231.

\bibitem[{Angulo et~al.(2016)Angulo, Ahmed, \protect\BIBand{}
  Dey}]{angulo2016improving}
Angulo G, Ahmed S, Dey SS (2016) {Improving the integer L-shaped method}.
  \emph{INFORMS Journal on Computing} 28(3):483--499.

\bibitem[{Benders(1962)}]{benders1962partitioning}
Benders JF (1962) Partitioning procedures for solving mixed-variables
  programming problems. \emph{Numerische Mathematik} 4(1):238--252.

\bibitem[{Bley et~al.(2013)Bley, D'Andreagiovanni, \protect\BIBand{}
  Karch}]{ble13}
Bley A, D'Andreagiovanni F, Karch D (2013) Scheduling technology migration in
  {WDM} networks. \emph{ITG Symposium on Photonic Networks}, 1--5.

\bibitem[{Bredstrom \protect\BIBand{}
  R{\"o}nnqvist(2007)}]{bredstrom2007abranch}
Bredstrom D, R{\"o}nnqvist M (2007) A branch and price algorithm for the
  combined vehicle routing and scheduling problem with synchronization
  constraints. Technical report, Department of Finance and Management Science,
  Norwegian School of Economics and Business Administration.

\bibitem[{Chvatal et~al.(1983)Chvatal, Chvatal et~al.}]{chvatal1983linear}
Chvatal V, Chvatal V, et~al. (1983) \emph{Linear programming} (Macmillan).

\bibitem[{Ciena(2013)}]{ciena2013}
Ciena (2013) The network modernization imperative.
  \urlprefix\url{https://media.ciena.com/documents/The-Network-Modernization-Imperative_page1.pdf}.

\bibitem[{Cir{\'e} et~al.(2016)Cir{\'e}, Coban, \protect\BIBand{}
  Hooker}]{cire2016logic}
Cir{\'e} AA, Coban E, Hooker JN (2016) Logic-based benders decomposition for
  planning and scheduling: A computational analysis. \emph{The Knowledge
  Engineering Review} 31(5):440--451.

\bibitem[{Dawadi et~al.(2021)Dawadi, Rawat, Joshi, Manzoni, \protect\BIBand{}
  Keitsch}]{dawadi2021migration}
Dawadi BR, Rawat DB, Joshi SR, Manzoni P, Keitsch MM (2021) Migration cost
  optimization for service provider legacy network migration to
  software-defined {IPv6} network. \emph{International Journal of Network
  Management} 31(4):e2145.

\bibitem[{Drexl(2012)}]{drexl2012synchronization}
Drexl M (2012) Synchronization in vehicle routing—a survey of {VRPs} with
  multiple synchronization constraints. \emph{Transportation Science}
  46(3):297--316.

\bibitem[{Eksioglu et~al.(2009)Eksioglu, Vural, \protect\BIBand{}
  Reisman}]{eksioglu2009vehicle}
Eksioglu B, Vural AV, Reisman A (2009) The vehicle routing problem: A taxonomic
  review. \emph{Computers \& Industrial Engineering} 57(4):1472--1483.

\bibitem[{Elci \protect\BIBand{} Hooker(2020)}]{elci2020stochastic}
Elci O, Hooker J (2020) Stochastic planning and scheduling with logic-based
  {Benders} decomposition. \emph{arXiv preprint arXiv:2012.14074} .

\bibitem[{Fawaz et~al.(2004)Fawaz, Daheb, Audouin, Du-Pond, \protect\BIBand{}
  Pujolle}]{fawaz2004service}
Fawaz W, Daheb B, Audouin O, Du-Pond M, Pujolle G (2004) Service level
  agreement and provisioning in optical networks. \emph{IEEE Communications
  Magazine} 42(1):36--43.

\bibitem[{{Gurobi Optimization, LLC}(2021)}]{gurobi}
{Gurobi Optimization, LLC} (2021) {Gurobi} optimizer reference manual.
  \urlprefix\url{https://www.gurobi.com}.

\bibitem[{Hashemi~Doulabi et~al.(2020)Hashemi~Doulabi, Pesant,
  \protect\BIBand{} Rousseau}]{hashemi2020vehicle}
Hashemi~Doulabi H, Pesant G, Rousseau LM (2020) Vehicle routing problems with
  synchronized visits and stochastic travel and service times: Applications in
  healthcare. \emph{Transportation Science} 54(4):1053--1072.

\bibitem[{Hojabri et~al.(2018)Hojabri, Gendreau, Potvin, \protect\BIBand{}
  Rousseau}]{hojabri2018large}
Hojabri H, Gendreau M, Potvin JY, Rousseau LM (2018) Large neighborhood search
  with constraint programming for a vehicle routing problem with
  synchronization constraints. \emph{Computers \& Operations Research}
  92:87--97.

\bibitem[{Hooker et~al.(2000)Hooker, Ottosson, Thorsteinsson, \protect\BIBand{}
  Kim}]{hooker2000scheme}
Hooker J, Ottosson G, Thorsteinsson ES, Kim HJ (2000) A scheme for unifying
  optimization and constraint satisfaction methods. \emph{The Knowledge
  Engineering Review} 15(1):11--30.

\bibitem[{Hooker \protect\BIBand{} Ottosson(2003)}]{hooker2003logic}
Hooker JN, Ottosson G (2003) {Logic-based Benders decomposition}.
  \emph{Mathematical Programming} 96(1):33--60.

\bibitem[{Hooker \protect\BIBand{} van Hoeve(2018)}]{hooker2018constraint}
Hooker JN, van Hoeve WJ (2018) Constraint programming and operations research.
  \emph{Constraints} 23(2):172--195.

\bibitem[{Hooker et~al.(2012)}]{hooker2012integrated}
Hooker JN, et~al. (2012) \emph{Integrated methods for optimization}, volume 170
  (Springer).

\bibitem[{{IBM® ILOG® CP® Optimizer}(2021)}]{cplex}
{IBM® ILOG® CP® Optimizer} (2021) {CPLEX CP Optimizer Users Manual}.
  \urlprefix\url{https://www.ibm.com/analytics/cplex-cp-optimizer}.

\bibitem[{Jain \protect\BIBand{} Grossmann(2001)}]{jain2001algorithms}
Jain V, Grossmann IE (2001) Algorithms for hybrid {MILP/CP} models for a class
  of optimization problems. \emph{INFORMS Journal on Computing} 13(4):258--276.

\bibitem[{Jaumard \protect\BIBand{} Pouya(2018)}]{jaumard2018migration}
Jaumard B, Pouya H (2018) Migration plan with minimum overall migration time or
  cost. \emph{Journal of Optical Communications and Networking} 10(1):1--13.

\bibitem[{Jaumard et~al.(2016)Jaumard, Pouya, Fahim, \protect\BIBand{}
  Barrios}]{brigitteHamed}
Jaumard B, Pouya H, Fahim R, Barrios A (2016) Planning network migration.
  \emph{IEEE International Conference on Communications - ICC}, 1--6.

\bibitem[{Javad-Kalbasi \protect\BIBand{} Valaee(2021)}]{javad2021efficient}
Javad-Kalbasi M, Valaee S (2021) Efficient migration to the next generation of
  networks based on {Digital Annealing}. \emph{ICASSP 2021-2021 IEEE
  International Conference on Acoustics, Speech and Signal Processing
  (ICASSP)}, 4740--4744 (IEEE).

\bibitem[{Knight et~al.(2011)Knight, Nguyen, Falkner, Bowden, \protect\BIBand{}
  Roughan}]{6027859}
Knight S, Nguyen H, Falkner N, Bowden R, Roughan M (2011) The internet topology
  zoo. \emph{Selected Areas in Communications, IEEE Journal on} 29(9):1765
  --1775, ISSN 0733-8716,
  \urlprefix\url{http://dx.doi.org/10.1109/JSAC.2011.111002}.

\bibitem[{Labadie et~al.(2014)Labadie, Prins, \protect\BIBand{}
  Yang}]{labadie2014iterated}
Labadie N, Prins C, Yang Y (2014) Iterated local search for a vehicle routing
  problem with synchronization constraints. \emph{International Conference on
  Operations Research and Enterprise Systems - ICORES}, 257--263.

\bibitem[{Laborie et~al.(2018)Laborie, Rogerie, Shaw, \protect\BIBand{}
  Vil{\'\i}m}]{laborie2018ibm}
Laborie P, Rogerie J, Shaw P, Vil{\'\i}m P (2018) {IBM} {ILOG} {CP Pptimizer}
  for scheduling. \emph{Constraints} 23(2):210--250.

\bibitem[{Laporte \protect\BIBand{} Louveaux(1993)}]{laporte1993integer}
Laporte G, Louveaux FV (1993) {The integer L-shaped method for stochastic
  integer programs with complete recourse}. \emph{Operations Research Letters}
  13(3):133--142.

\bibitem[{Li et~al.(2020)Li, Qin, Baldacci, \protect\BIBand{}
  Zhu}]{li2020branch}
Li J, Qin H, Baldacci R, Zhu W (2020) Branch-and-price-and-cut for the
  synchronized vehicle routing problem with split delivery, proportional
  service time and multiple time windows. \emph{Transportation Research Part E:
  Logistics and Transportation Review} 140:101955.

\bibitem[{Loken et~al.(2010)Loken, Gruner, Groer, Peltier, Bunn, Craig,
  Henriques, Dempsey, Yu, Chen et~al.}]{loken2010scinet}
Loken C, Gruner D, Groer L, Peltier R, Bunn N, Craig M, Henriques T, Dempsey J,
  Yu CH, Chen J, et~al. (2010) {SciNet: lessons learned from building a
  power-efficient top-20 system and data centre}. \emph{Journal of
  Physics-Conference Series}, volume 256, 012026.

\bibitem[{Matousek \protect\BIBand{}
  G{\"a}rtner(2007)}]{matousek2007understanding}
Matousek J, G{\"a}rtner B (2007) \emph{Understanding and using linear
  programming} (Springer Science \& Business Media).

\bibitem[{Perron \protect\BIBand{} Furnon(2019)}]{ortools}
Perron L, Furnon V (2019) {OR-Tools}.
  \urlprefix\url{https://developers.google.com/optimization/}.

\bibitem[{Podhradsky(2004)}]{podhradsky2004}
Podhradsky P (2004) Migration scenarios and convergence processes towards {NGN}
  (present state and future trends). \emph{IEEE International Symposium
  Electronics in Marine}, 39--46.

\bibitem[{Ponce et~al.(2019)Ponce, van Zon, Northrup, Gruner, Chen, Ertinaz,
  Fedoseev, Groer, Mao, Mundim et~al.}]{ponce2019deploying}
Ponce M, van Zon R, Northrup S, Gruner D, Chen J, Ertinaz F, Fedoseev A, Groer
  L, Mao F, Mundim BC, et~al. (2019) {Deploying a top-100 supercomputer for
  large parallel workloads: The Niagara supercomputer}. \emph{Proceedings of
  the Practice and Experience in Advanced Research Computing on Rise of the
  Machines (learning)}, 1--8 (Association for Computing Machinery, New York,
  NY, United States).

\bibitem[{Poularakis et~al.(2019)Poularakis, Iosifidis, Smaragdakis,
  \protect\BIBand{} Tassiulas}]{poularakis2019optimizing}
Poularakis K, Iosifidis G, Smaragdakis G, Tassiulas L (2019) Optimizing gradual
  {SDN} upgrades in {ISP} networks. \emph{IEEE/ACM Transactions on Networking}
  27(1):288--301.

\bibitem[{Pouya(2018)}]{pouya2018new}
Pouya H (2018) \emph{New Models and Algorithms in Telecommunication Networks}.
  Ph.D. thesis, Concordia University.

\bibitem[{Pouya \protect\BIBand{} Jaumard(2017)}]{pouya2017efficient}
Pouya H, Jaumard B (2017) Efficient network migration planning. \emph{2017 19th
  International Conference on Transparent Optical Networks (ICTON)}, 1--4
  (IEEE).

\bibitem[{Pouya et~al.(2017)Pouya, Jaumard, \protect\BIBand{}
  Preston-Thomas}]{pouya2017minimum}
Pouya H, Jaumard B, Preston-Thomas C (2017) Minimum network migration cost and
  duration. \emph{IEEE Sarnoff Symposium}, 1--6.

\bibitem[{Reinhardt et~al.(2013)Reinhardt, Clausen, \protect\BIBand{}
  Pisinger}]{reinhardt2013synchronized}
Reinhardt LB, Clausen T, Pisinger D (2013) Synchronized dial-a-ride
  transportation of disabled passengers at airports. \emph{European Journal of
  Operational Research} 225(1):106--117.

\bibitem[{Roshanaei et~al.(2017)Roshanaei, Luong, Aleman, \protect\BIBand{}
  Urbach}]{roshanaei2017propagating}
Roshanaei V, Luong C, Aleman DM, Urbach D (2017) Propagating logic-based
  {Benders} decomposition approaches for distributed operating room scheduling.
  \emph{European Journal of Operational Research} 257(2):439--455.

\bibitem[{Rousseau et~al.(2004)Rousseau, Gendreau, Pesant, \protect\BIBand{}
  Focacci}]{rousseau2004solving}
Rousseau LM, Gendreau M, Pesant G, Focacci F (2004) Solving {VRPTWs} with
  constraint programming based column generation. \emph{Annals of Operations
  Research} 130(1-4):199--216.

\bibitem[{Salazar-Aguilar et~al.(2013)Salazar-Aguilar, Langevin,
  \protect\BIBand{} Laporte}]{salazar2013synchronized}
Salazar-Aguilar MA, Langevin A, Laporte G (2013) The synchronized arc and node
  routing problem: Application to road marking. \emph{Computers \& Operations
  Research} 40(7):1708--1715.

\bibitem[{T{\"u}rk et~al.(2012)T{\"u}rk, Liu, Radeke, \protect\BIBand{}
  Lehnert}]{turk2012networkMigOpt}
T{\"u}rk S, Liu Y, Radeke R, Lehnert R (2012) Network migration optimization
  using genetic algorithms. \emph{Information and Communication Technologies},
  volume 7479, 112--123.

\bibitem[{Vance(1998)}]{vance1998branch}
Vance PH (1998) Branch-and-price algorithms for the one-dimensional cutting
  stock problem. \emph{Computational Optimization and Applications}
  9(3):211--228.

\bibitem[{Vanderbeck(2011)}]{vanderbeck2011branching}
Vanderbeck F (2011) Branching in branch-and-price: a generic scheme.
  \emph{Mathematical Programming} 130(2):249--294.

\end{thebibliography}

\newpage

\ECSwitch

%\ECDisclaimer
%%%%%%%%%%%%%%%%%%%%%%%%%%%%%%%%%%%%%%%%%%%%%%%%%%%%%%%%%%

%%% Main head for the e-companion
\ECHead{e-Companion}
\section{Auxiliary MIP Pricing Problem}
In addition to $n_{ss'}$ and $n_{\textsc{cir}}$ described in Section \ref{sec:opt_models}, we define the new decision variables in Table \ref{tab:CG_Sub_var}.
Denote by $\overline{\boldsymbol{\dual}}^{\eqref{eq: m_shift_eq_NMP_CG_RMP}}, \overline{\boldsymbol{\dual}}^{\eqref{eq: techs_in_r_mw_NMP_CG_RMP}},  \boldsymbol{\dual}^{\eqref{eq: engineers_NMP_CG_RMP}}$, the optimal dual solutions associated with constraints ${\eqref{eq: m_shift_eq_NMP_CG_RMP}}$, ${\eqref{eq: techs_in_r_mw_NMP_CG_RMP}}$, ${\eqref{eq: engineers_NMP_CG_RMP}} $, respectively. The pricing problem generating a shift for $ (r,w) $ is:\\[-9mm]
\begin{subequations}
\begin{align}
\text{SP}^{'\text{CG}}_{rw} = \min \ \  &\cost\durshift -\sum_{s\in \siteSet_r}\sum_{s'\in \siteSet} n_{ss'}\overline{\dual}_{ss'}^{\eqref{eq: m_shift_eq_NMP_CG_RMP}} - \overline{\dual}^{\eqref{eq: techs_in_r_mw_NMP_CG_RMP}}_{rw} - \overline{\dual}^{\eqref{eq: engineers_NMP_CG_RMP}}_{w} \label{eq:pp-obj-app}\\
\text{s.t.}\ \  &h_{ss'} + h_{s's} \leq 1 && s,s'\in \siteSet_r, s\neq s' \label{eq: each_tech_one_ep_pricing_II}\\
&h_s \leq \sum_{s' \in S} h_{ss'} \leq M h_{s} \qquad&& s\in \siteSet_r \label{eq: sites_pricing_II}\\
&h_{ss'} \leq n_{ss'} \leq Mh_{ss'} \qquad&& s \in \siteSet_r, s' \in S, s\neq s' \label{eq: migrated_cir_pricing_II}\\
& \sum_{s\in \siteSet_r} t_{s_{\SRC},s} = 1, \quad \sum_{s\in \siteSet_r} t_{s,s_{\DST}} = 1 && \label{eq: exactly-one-site-for-end_pricing_I}\\
&t_{ss'} \leq \frac{1}{2}(h_s+h_{s'}) &&  s, s' \in \siteSet_{r}, s'>s 	\label{eq: t-s-sprime_pricing_I}\\
& \sum_{\substack{s' \in \siteSet_r^+\\ s' > s}} t_{ss'} + \sum_{\substack{s' \in \siteSet_r^+\\ s' < s}} t_{s's}  = 2 h_s \quad &&  s \in \siteSet_r \label{eq: at-most-two-site-after-and-before_pricing_I}\\
& \sum_{s\in S_r} \sum_{\substack{s' \in \siteSet_r\\ s' > s}} t_{ss'} = \sum\limits_{s\in \siteSet_r} h_s - 1 \label{eq:path_pricing_I}\\
&\durmigr \sum\limits_{s\in S_r} \sum_{\substack{s'\in S \\ s\neq s'}} n_{ss'} + \sum\limits_{s\in S_r} \sum\limits_{\substack{s'\in S_r \\ s' > s}} \travel_{ss'} t_{ss'}\leq \durshift&&  \label{eq: dur_pricing_II}\\
&\sum_{\delta \in \duration} x_{\delta} = 1 \label{eq: one_dur_pricing_I}\\
&\sum_{\delta \in \duration} \duration_{\delta}x_{\delta} = \durshift \label{eq: final_dur_pricing_I}\\
& n_{ss'}, h_{ss'} \in \mathbb{Z}_+, h_s \in \{0,1\} \qquad&& s \in \siteSet_r, s' \in \siteSet, s\neq s' 
\end{align}
\end{subequations}
where $ M\geq 0 $ is a large number and $\siteSet^+ = \siteSet \cup \{ s_{\src}, s_{\dst} \}$, $\siteSet_r^+ = S_r \cup \{ s_{\src}, s_{\dst} \}$. 
\begin{table}[bth]
\centering
\caption{New decision variables for the CG subproblem}
\label{tab:CG_Sub_var}
  {\renewcommand{\arraystretch}{1.1}
  \small
\begin{tabular}{llp{5in}}
\toprule
\textbf{Variable} & \textbf{\small Type} & \textbf{\small Description}\\
\midrule
	$ h_{ss'} $ & Binary & $=1$ if the technician migrates at least one circuit endpoint in site $ s $ with the other endpoint in site $ s' $.\\
	$ h_s $ & Binary & $=1$ if the technician works in site $ s $ in this shift, $ 0 $ otherwise.\\
	$t_{ss'} $ & Binary & $=1$ if a travel from site $s$ to site $s'$ occurs in the shift under construction, $0$ otherwise.\\
	$ x_\delta $ & Binary & $=1$ if the length of the shift under construction is equal to $ \duration_\delta $, $0$ otherwise.\\
\bottomrule
\end{tabular}
}
\end{table}
$s_{\src}$ and $s_{\dst}$ are two dummy sites introduced as the beginning and end of a path. Objective function \eqref{eq:pp-obj-app} is the reduced-cost. 
Constraints \eqref{eq: each_tech_one_ep_pricing_II} ensure that at most one of the endpoints of every circuit $ c\in \cirSet_{ss'} $ can be migrated in every shift. Constraints \eqref{eq: sites_pricing_II} determine the sites where a technician works in the current shift. Constraints \eqref{eq: migrated_cir_pricing_II} assure that all migrated circuit endpoints are from the sites where the technician works in the current shift. Thanks to constraints \eqref{eq: exactly-one-site-for-end_pricing_I}, generated path in this configuration starts form dummy site  $s_{\src}$ and ends in dummy site $s_{\dst}$.
Constraints \eqref{eq: t-s-sprime_pricing_I} guarantee that if a travel occurs between the two sites $s$ and $s'$, the technician works in both.
Constraints \eqref{eq: at-most-two-site-after-and-before_pricing_I} determine the sites visited before and after every site $s\in \siteSet_r$.
Constraint \eqref{eq:path_pricing_I} ensures that we have a path linking all visited sites. Constraint \eqref{eq: dur_pricing_II} makes sure that the shift duration does not exceed the maximum predefined duration. Constraints \eqref{eq: one_dur_pricing_I} and \eqref{eq: final_dur_pricing_I} together determine the duration of the shift. The rest of the constraints determine the variables domain.

\section{Detailed Numerical Results}
\begin{table}[htbp]
  \centering
  \caption{Characteristics of instances}
  \scalebox{0.85}{\renewcommand{\arraystretch}{1.25}
  \small
\begin{tabular}{lrcrcrcrcrcrc}
\toprule
\multirow{3}[6]{*}{Dataset} & \multicolumn{12}{c}{Instance} \\
\cmidrule{2-13}      & \multicolumn{2}{c}{1} & \multicolumn{2}{c}{2} & \multicolumn{2}{c}{3} & \multicolumn{2}{c}{4} & \multicolumn{2}{c}{5} & \multicolumn{2}{c}{6} \\
\cmidrule(lr){2-3}\cmidrule(lr){4-5}\cmidrule(lr){6-7}\cmidrule(lr){8-9}\cmidrule(lr){10-11}\cmidrule(lr){12-13}      & \multicolumn{1}{l}{Endpoints} & $|\sitePair|$ & \multicolumn{1}{l}{Endpoints} & $|\sitePair|$ & \multicolumn{1}{l}{Endpoints} & $|\sitePair|$ & \multicolumn{1}{l}{Endpoints} & $|\sitePair|$ & \multicolumn{1}{l}{Endpoints} & $|\sitePair|$ & \multicolumn{1}{l}{Endpoints} & $|\sitePair|$ \\
\midrule
\EUnet5 & 66    & 27    & 60    & 26    & 46    & 22    & 66    & 26    & 66    & 30    & 66    & 32 \\
\EUnet6 & 70    & 30    & 72    & 32    & 84    & 31    & 88    & 34    & 80    & 33    & 74    & 28 \\
\EUnet7 & 98    & 40    & 84    & 36    & 94    & 34    & 80    & 35    & 90    & 39    & 86    & 35 \\
\EUnet8 & 116   & 41    & 94    & 40    & 110   & 43    & 102   & 38    & 106   & 42    & 108   & 39 \\
\EUnet9 & 110   & 46    & 150   & 44    & 112   & 45    & 112   & 41    & 104   & 42    & 104   & 42 \\
\midrule
\NextGen5 & 54    & 26    & 52    & 23    & 56    & 26    & 74    & 31    & 66    & 31    & 56    & 26 \\
\NextGen6 & 86    & 39    & 62    & 26    & 88    & 34    & 80    & 34    & 96    & 39    & 74    & 35 \\
\NextGen7 & 82    & 36    & 86    & 35    & 90    & 35    & 92    & 37    & 92    & 35    & 78    & 31 \\
\NextGen8 & 122   & 46    & 126   & 47    & 108   & 44    & 92    & 35    & 118   & 44    & 98    & 40 \\
\NextGen9 & 134   & 45    & 126   & 51    & 106   & 41    & 124   & 45    & 114   & 45    & 98    & 44 \\
\midrule
\Pionier5 & 78    & 38    & 82    & 39    & 70    & 35    & 90    & 39    & 98    & 43    & 74    & 36 \\
\Pionier6 & 112   & 48    & 112   & 48    & 106   & 51    & 106   & 47    & 96    & 40    & 104   & 42 \\
\midrule
\Sago5 & 76    & 35    & 64    & 29    & 88    & 39    & 74    & 33    & 82    & 35    & 90    & 39 \\
\Sago6 & 86    & 37    & 94    & 40    & 94    & 40    & 80    & 32    & 86    & 38    & 90    & 40 \\
\midrule
\Savvis5 & 72    & 35    & 98    & 42    & 102   & 41    & 80    & 36    & 54    & 25    & 86    & 40 \\
\Savvis6 & 104   & 40    & 90    & 40    & 84    & 40    & 106   & 41    & 88    & 38    & 94    & 40 \\
\Savvis7 & 120   & 46    & 106   & 44    & 88    & 37    & 114   & 47    & 114   & 45    & 92    & 39 \\
\midrule
\VisionNet5 & 86    & 40    & 100   & 47    & 104   & 45    & 88    & 38    & 102   & 44    & 82    & 36 \\
\VisionNet6 & 116   & 52    & 106   & 51    & 120   & 52    & 120   & 47    & 126   & 52    & 110   & 49 \\
\bottomrule
\end{tabular}%
    }
  \label{tab:instances-sitepairs}%
\end{table}%

\begin{table}[htbp]
  \centering
  \caption{Breakdown of solution times}
  \scalebox{0.85}{\renewcommand{\arraystretch}{1.25}
  \small
    \begin{tabular}{lrcrrrrrrrr}
    \toprule
    \multirow{2}[2]{*}{Dataset} & \multicolumn{2}{c}{Iterations} & \multicolumn{2}{c}{Sol Time} & \multicolumn{2}{c}{LBBD MP Time} & \multicolumn{2}{c}{CG Time} & \multicolumn{2}{c}{CP Time} \\
    \cmidrule(lr){2-3}\cmidrule(lr){4-5}\cmidrule(lr){6-7}\cmidrule(lr){8-9}\cmidrule(lr){10-11}
          & Mean  & CI width & Mean  & CI width & Mean  & CI width & Mean  & CI width & Mean  & CI width \\
    \midrule
    \EUnet5 & 3.5   & 0.6   & 1124.1 & 268.6 & 548.2 & 198.3 & 540.1 & 196.0 & 575.8 & 147.6 \\
    \EUnet6 & 10.5  & 4.3   & 4607.1 & 2309.7 & 2953.3 & 1929.8 & 2939.2 & 1928.9 & 1653.5 & 680.5 \\
    \EUnet7 & 6.5   & 2.0   & 1459.0 & 509.7 & 45.7  & 23.4  & 27.6  & 13.7  & 1413.2 & 497.5 \\
    \EUnet8 & 7.3   & 1.6   & 2582.6 & 539.9 & 693.8 & 271.9 & 531.6 & 225.2 & 1888.4 & 563.5 \\
    \EUnet9 & 10.8  & 6.1   & 5109.8 & 3376.7 & 2201.9 & 1889.0 & 1384.5 & 1000.8 & 2907.4 & 1817.8 \\
    \midrule
    \NextGen5 & 6.0   & 3.0   & 699.2 & 560.3 & 18.0  & 16.1  & 4.6   & 2.7   & 680.9 & 561.0 \\
    \NextGen6 & 4.8   & 2.3   & 675.7 & 374.2 & 94.8  & 42.8  & 17.6  & 6.3   & 580.8 & 361.1 \\
    \NextGen7 & 5.3   & 2.0   & 997.9 & 415.8 & 177.6 & 63.9  & 32.9  & 15.0  & 820.1 & 381.1 \\
    \NextGen8 & 8.2   & 4.1   & 2445.6 & 399.6 & 1182.4 & 626.9 & 135.9 & 66.9  & 1262.8 & 433.9 \\
    \NextGen9 & 8.0   & 2.9   & 3120.4 & 1939.1 & 1980.7 & 1546.1 & 67.1  & 42.6  & 1139.5 & 543.4 \\
    \midrule
    \Pionier5 & 5.7   & 1.6   & 1817.2 & 796.4 & 812.9 & 505.1 & 726.2 & 475.5 & 1004.1 & 322.2 \\
    \Pionier6 & 9.2   & 2.0   & 4929.4 & 1385.8 & 3267.2 & 1306.1 & 2961.9 & 1167.6 & 1661.8 & 272.9 \\
    \midrule
    \Sago5 & 8.3   & 2.3   & 2387.6 & 697.1 & 836.1 & 293.6 & 748.0 & 319.1 & 1551.2 & 480.3 \\
    \Sago6 & 5.7   & 2.3   & 3635.1 & 2640.4 & 2630.0 & 2194.2 & 2473.5 & 2132.5 & 1005.0 & 490.1 \\
    \midrule
    \Savvis5 & 4.5   & 1.4   & 619.6 & 475.2 & 400.5 & 408.2 & 27.1  & 12.1  & 218.8 & 94.4 \\
    \Savvis6 & 6.0   & 1.6   & 855.6 & 460.3 & 608.2 & 393.2 & 36.1  & 15.3  & 247.1 & 108.7 \\
    \Savvis7 & 6.0   & 3.7   & 2472.6 & 2731.3 & 2184.8 & 2503.3 & 96.2  & 108.1 & 287.3 & 251.1 \\
    \midrule
    \VisionNet5 & 5.8   & 3.5   & 2259.2 & 899.9 & 1217.2 & 523.3 & 1112.4 & 507.2 & 1041.8 & 541.1 \\
    \VisionNet6 & 7.3   & 3.1   & 3308.0 & 769.4 & 1759.4 & 718.7 & 139.2 & 36.3  & 1548.3 & 996.6 \\
    \bottomrule
    \end{tabular}%
    }
  \label{tab:performance}%
\end{table}%

\begin{table}[htbp]
  \centering
  \caption{The number of cuts and columns generated during the solution process}
    \scalebox{0.85}{\renewcommand{\arraystretch}{1.25}
  \small
    \begin{tabular}{lrrrrrrrr}
    \toprule
    \multirow{2}[4]{*}{Dataset} & \multicolumn{2}{c}{\#Cut$^{\textsc{LBBD}}_{\textsc{opt}}$} & \multicolumn{2}{c}{\#Cut$^{\textsc{LBBD}}_{\textsc{feas}}$} & \multicolumn{2}{c}{\#Cut$^{\textsc{BD}}_{\textsc{opt}}$} & \multicolumn{2}{c}{\#Columns} \\
\cmidrule(lr){2-3} \cmidrule(lr){4-5}\cmidrule(lr){6-7} \cmidrule(lr){8-9}           & Mean  & CI width & Mean  & CI width & Mean  & CI width & Mean  & CI width \\
    \midrule
    \EUnet5 & 121.8 & 35.3  & 27.3  & 11.2  & 33.0  & 11.6  & 408.8 & 180.2 \\
    \EUnet6 & 566.8 & 273.0 & 84.0  & 30.5  & 79.2  & 43.0  & 951.2 & 375.3 \\
    \EUnet7 & 361.2 & 101.1 & 66.3  & 44.1  & 76.8  & 36.8  & 1350.8 & 565.1 \\
    \EUnet8 & 463.7 & 114.1 & 88.3  & 29.0  & 170.8 & 62.0  & 3693.8 & 1120.1 \\
    \EUnet9 & 919.2 & 605.7 & 119.3 & 104.5 & 480.7 & 464.9 & 7318.7 & 4824.6 \\
    \midrule
    \NextGen5 & 218.3 & 136.2 & 59.2  & 23.6  & 35.3  & 28.4  & 45.5  & 0.9 \\
    \NextGen6 & 237.5 & 167.3 & 42.7  & 25.4  & 115.8 & 43.8  & 159.8 & 162.7 \\
    \NextGen7 & 267.7 & 130.6 & 61.5  & 40.4  & 220.5 & 50.3  & 274.8 & 411.2 \\
    \NextGen8 & 434.3 & 148.6 & 66.0  & 31.4  & 386.2 & 164.1 & 1405.5 & 646.1 \\
    \NextGen9 & 622.3 & 317.5 & 103.7 & 32.8  & 666.5 & 405.6 & 1456.2 & 991.8 \\
    \midrule
    \Pionier5 & 289.0 & 96.7  & 76.0  & 20.5  & 110.0 & 73.0  & 741.5 & 462.9 \\
    \Pionier6 & 666.8 & 159.0 & 124.8 & 36.9  & 308.2 & 70.9  & 2101.8 & 1108.5 \\
    \midrule
    \Sago5 & 476.8 & 185.9 & 78.3  & 36.1  & 232.7 & 165.6 & 1026.8 & 691.0 \\
    \Sago6 & 307.7 & 161.8 & 37.0  & 17.9  & 365.5 & 260.5 & 1024.3 & 582.1 \\
    \midrule
    \Savvis5 & 215.8 & 83.0  & 50.0  & 29.9  & 176.0 & 107.9 & 57.0  & 0.0 \\
    \Savvis6 & 343.3 & 137.5 & 72.2  & 13.6  & 219.0 & 83.0  & 57.0  & 0.0 \\
    \Savvis7 & 375.7 & 320.2 & 65.2  & 42.6  & 661.0 & 728.0 & 57.0  & 0.0 \\
    \midrule
    \VisionNet5 & 385.5 & 322.0 & 63.3  & 43.9  & 136.7 & 53.6  & 1772.0 & 1091.4 \\
    \VisionNet6 & 571.3 & 328.8 & 115.8 & 53.0  & 411.5 & 102.8 & 4607.8 & 1571.0 \\
    \bottomrule
    \end{tabular}%
    }
  \label{tab:number-cuts}%
\end{table}%

\begin{figure}[htbp]
	\centering 
	\subfloat[\# Shifts vs. cost\label{fig:shift_cost_sago}]{\scalebox{0.57}{
	\input{Figures/Sago_ShiftCost}}}
	% -------------------
	\subfloat[Distribution of shifts\label{fig:shift_dur_sago}]{\scalebox{0.57}{
	\input{Figures/EU_Shifts}
	}} \\
    % -------------------
	\subfloat[Cost vs. efficient working time\label{fig:working_cost_sago}]{\scalebox{0.6}{
	\input{Figures/Sago_WorkingCost}
	}}\\
	% -------------------
	\caption{Cost analysis for \Sago}
	\label{fig:cost_analysis_sago}
\end{figure}

\begin{figure}[htbp]
	\centering 
	\subfloat[\# Shifts vs. cost\label{fig:shift_cost_savvis}]{\scalebox{0.57}{
	\input{Figures/Saavi_ShiftCost}}}
	% -------------------
	\subfloat[Distribution of shifts\label{fig:shift_dur_savvis}]{\scalebox{0.57}{
	\input{Figures/Saavi_Shifts}
	}} \\
    % -------------------
	\subfloat[Cost vs. efficient working time\label{fig:working_cost_savvis}]{\scalebox{0.6}{
	\input{Figures/Saavi_WorkingCost}
	}}\\
% 	-------------------
	\caption{Cost analysis for \Savvis}
	\label{fig:cost_analysis_savvis}
\end{figure}

\end{document}